\documentclass[a4 paper, 12pt]{article} 
\usepackage{amssymb,amsthm}
\usepackage{amsmath}
\usepackage{subcaption}
\DeclareMathOperator*{\argmin}{arg\,min}

\usepackage{tikz}
\usepackage{graphicx}
\graphicspath{{spindlefigvis1/}}
\usetikzlibrary{decorations.pathmorphing}
\tikzset{snake it/.style={decorate, decoration=snake}}

\newcommand{\projecttitle}{Three dimensional Compton scattering tomography}
\newcommand{\projectauthor}{By James Webber and William Lionheart}
\theoremstyle{plain}
\newtheorem{theorem}{Theorem}
\newtheorem{lemma}{Lemma}
\theoremstyle{definition}

\newtheorem{corollary}{Corollary}

\def\titlep{\thispagestyle{empty}\null\vskip1in\begin{center}
           \Huge\uppercase{\projecttitle}
\end{center}
\vspace{10mm}
\begin{center}
\text{\projectauthor}
\end{center}
\vspace{10mm}
}

\begin{document}
\titlep
\begin{abstract}
We propose a new acquisition geometry for electron density reconstruction in three dimensional X-ray Compton imaging using a monochromatic source. This leads us to a new three dimensional inverse problem where we aim to reconstruct a real valued function $f$ (the electron density) from its integrals over spindle tori. We prove injectivity of a generalized spindle torus transform on the set of smooth functions compactly supported on a hollow ball. This is obtained through the explicit inversion of a class of Volterra integral operators, whose solutions give us an expression for the harmonic coefficients of $f$. The polychromatic source case is later considered, and we prove injectivity of a new spindle interior transform, apple transform and apple interior transform on the set of smooth functions compactly supported on a hollow ball. 

A possible physical model is suggested for both source types. We also provide simulated density reconstructions with varying levels of added pseudo random noise and model the systematic error due to the attenuation of the incoming and scattered rays in our simulation.
\end{abstract}

\section{Introduction}
In this paper we lay the foundations for a new three dimensional imaging technique in X-ray Compton scattering tomography. Recent publications present various two dimensional scattering modalities, where a function in the plane is reconstructed from its integrals over circular arcs \cite{pal1,NT,norton}. Three dimensional Compton tomography is also considered in the literature, where a gamma source is reconstructed from its integrals over cones with a fixed axis direction \cite{cone,cone1,cone2}. In \cite{CST}, Truong and Nguyen give a history of Compton scattering tomography, from the point by point reconstruction case in earlier modalities to the circular arc transform modalities in later work. Here we present a new three dimensional scattering modality, where we aim to reconstruct the electron density (the number of electrons per unit volume) from its integrals over the surfaces of revolution of circular arcs. This work provides the theoretical basis for a new form of non invasive density determination which would be applied in fields such as fossil imaging, airport baggage screening and more generally in X-ray spectroscopic imaging. Our main goal is to show that a unique three dimensional density reconstruction is possible with knowledge of the Compton scattered intensity with our proposed acquisition geometry, and to provide an analytic expression for the density in terms of the Compton scattered data.

Compton scattering is the process which describes the inelastic scattering of photons with charged particles (usually electrons). A loss in photon energy occurs upon the collision. This is known as the Compton effect. The energy loss is dependant on the initial photon energy and the angle of scattering and is described by the following equation:
\begin{equation}
\label{equ1}
E_s=\frac{E_{\lambda}}{1+\left(E_{\lambda}/E_0\right)\left(1-\cos\omega\right)}.
\end{equation}
Here $E_s$ is the energy of the scattered photon which had an initial energy $E_{\lambda}$, $\omega$ is the scattering angle and $E_0\approx 511$keV is the electron rest energy. Typically Compton scattering refers to the scattering of photons in the mid energy range. That is, the scattering of X-rays and gamma rays with photon energies ranging from 1keV up to 1MeV. Forward Compton scattering is the scattering of photons with scattering angles $\omega\leq \pi/2$. Conversely, backscatter refers to the scattering of photons at angles $\omega>\pi/2$.

\begin{figure}[!h]
\centering
\begin{tikzpicture}[scale=4]
\path ({-cos(30)},0.5) coordinate (S);
\path ({cos(30)}, 0.5) coordinate (D);
\path (-0.5,0.86) coordinate (w);
\path (-0.13,1.22) coordinate (a);
\draw [domain=30:150] plot ({cos(\x)}, {sin(\x)});
\draw [thin, dashed] (S) -- (D);
\draw (S) -- (w);
\draw (w) -- (D);
\draw [thin,dashed] (w)--(a);
\node at (-0.95,0.48) {$s$};
\node at (0.95,0.5) {$d$};
\node at (-0.35,0.9) {$\omega$};
\node at (0.9,0.9) {$C$};
\end{tikzpicture}
\caption{A circular arc $C$ is the locus of scattering points for a given measured energy $E_s$ for source and detector positions $s$ and $d$.}
\label{fig1}
\end{figure}
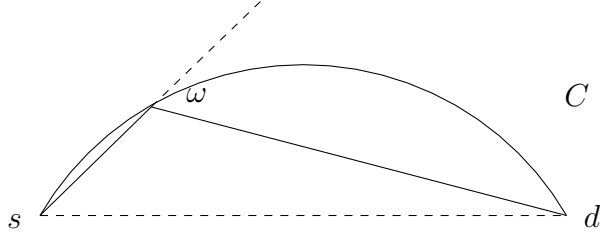

If the photon source is monochromatic ($E_{\lambda}$ is fixed) and the detector is energy resolving (we can measure photon intensity at a given scattered energy $E_s$), then, in a given plane, the locus of scattering points is a circular arc connecting the source and detector points. See figure \ref{fig1}.

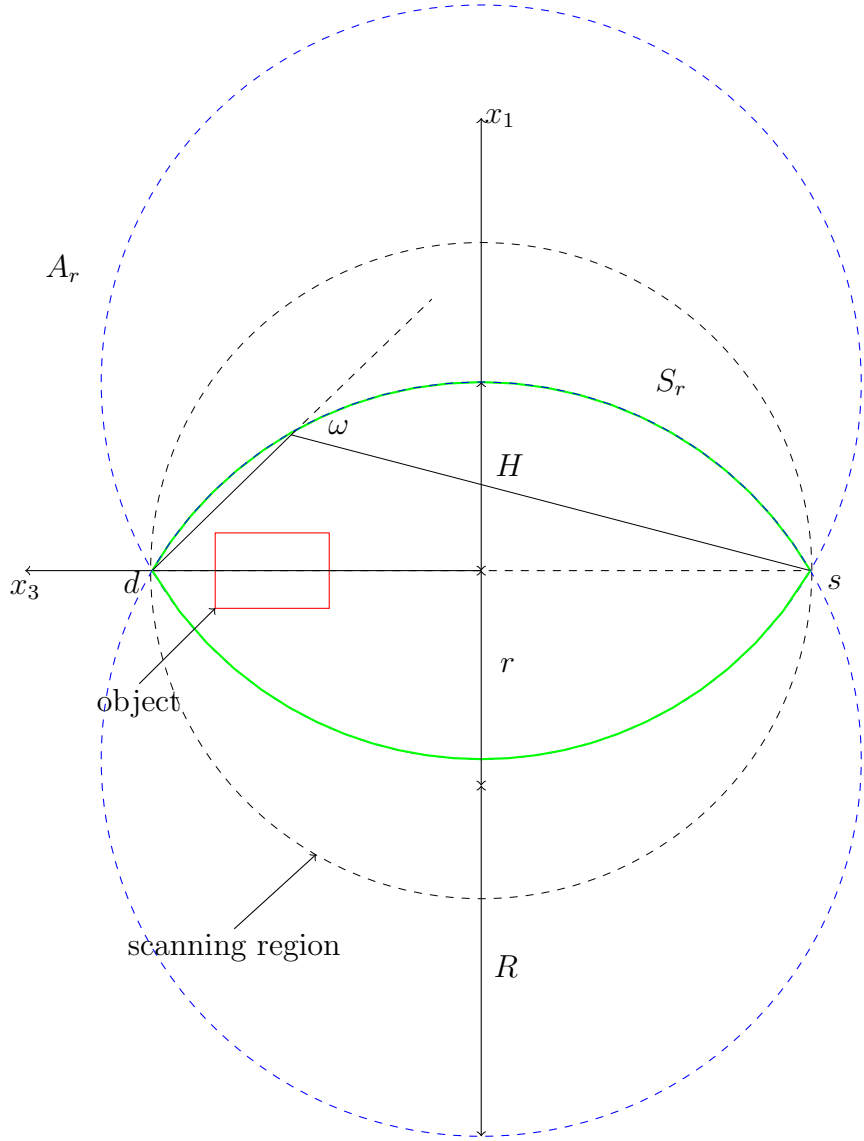
\begin{figure}[!h]
\centering
\begin{tikzpicture}[scale=5]
\path ({-cos(30)},0.5) coordinate (S);
\path ({cos(30)}, 0.5) coordinate (D);
\path (-0.5,0.86) coordinate (w);
\path (-0.13,1.22) coordinate (a);
\draw [green, thick, domain=30:150] plot ({cos(\x)}, {sin(\x)});
\draw [thin, dashed] (S) -- (D);
\draw [very thin] (S) -- (w);
\draw [very thin] (w) -- (D);
\draw [thin,dashed] (w)--(a);
\node at (-0.92,0.47) {$d$};
\node at (0.93,0.47) {$s$};
\node at (-0.375,0.88) {$\omega$};
\draw [very thin, dashed] (0,0.5) circle [radius=0.87];
\draw [red] (-0.7,0.4) rectangle (-0.4,0.6);
\draw [<->] (0,1)--(0,0.5);
\node at (0.075,0.78) {$H$};
\draw [->] (0,1)--(0,1.7);
\draw [->] (0,0.5)--(-1.2,0.5);
\node at (0.05,1.7) {$x_1$};
\node at (-1.2,0.45) {$x_3$};
\node at (0.5,1) {$S_r$};
\node at (-1.1,1.3) {$A_r$};
\draw [blue,dashed] (0,0) circle [radius=1];
\draw [blue,dashed] (0,1) circle [radius=1];
\draw [green,thick,domain=-0.866:0.866] plot(\x, {-sqrt(1-pow(\x,2))+1});
\draw [->] (-0.9, 0.2)->(-0.7,0.4);
\node at (-0.9, 0.15) {object};
\draw [->] (-0.65,-0.45)->(-0.435,-0.253);
\node at (-0.65,-0.5) {scanning region};
\draw [<->] (0,0.5)--(0,-0.07);
\node at (0.07,0.25) {$r$};
\draw [<->] (0,-0.07)--(0,-1);
\node at (0.065,-0.55) {$R$};
\end{tikzpicture}
\caption{A spindle torus slice with height $H$, tube radius $R$ and tube centre offset $r$ has tips at source and detector points $s$ and $d$. The spindle $S$ and the apple $A$ are highlighted by a solid green line and a dashed blue line respectively. The spherical scanning region highlighted has unit radius.}
\label{fig2}
\end{figure}

A \textit{spindle torus} is the surface of revolution of a circular arc. Specifically we define:
\begin{equation}
T_{r}=\left\{(x_1,x_2,x_3)\in \mathbb{R}^3 : \left(r-\sqrt{x_1^2+x_2^2}\right)^2+x_3^2=1+r^2\right\}
\end{equation}
to be the spindle torus, radially symmetric about the $x_3$ axis, with tube centre offset $r\geq 0$ and tube radius $\sqrt{1+r^2}$. Let $B_{0,1}$ denote the unit ball in $\mathbb{R}^3$. Then we define the \textit{spindle} $S_r=T_r \cap B_{0,1}$ as the interior of the torus $T_r$ and we define the \textit{apple} $A_r=T_r\backslash S_r$ as the remaining exterior. See figure \ref{fig2}.

In three dimensions, the surface of scatterers is the surface of revolution of a circular arc $C$ about its circle chord $sd$. Equivalently, the surface of scattering points is a spindle torus with an axis of revolution $sd$, whose tube radius and tube centre offset are determined by the distance $|sd|$ and the scattering angle $\omega$. In \cite{NT}, Nguyen and Truong present an acquisition geometry in two dimensions for a monochromatic source (e.g. a gamma ray source) and energy resolving detector pair, where the source and detector remain at a fixed distance opposite one another and are rotated about the origin on the curve $S^1$ (the unit circle). Here, the dimensionality of the data is two (an energy variable and a one dimensional rotation). Taking our inspiration from Nguyen and Truong's idea, we propose a novel acquisition geometry in three dimensions for a single source and energy resolving detector pair, which are rotated opposite each other at a fixed distance about the origin on the surface $S^2$ (the unit sphere). Our data set is three dimensional (an energy variable and a two dimensional rotation). 

As illustrated in figure \ref{fig2}, the forward scattered (for scattering angles $\omega\leq \pi/2$) intensity measured at the detector $d$ for a given energy $E_s$ can be written as a weighted integral over the spindle $S_r$ (the measured energy $E_s$ determines $r$). With this in mind we aim to reconstruct a function  supported within the unit ball from its weighted integrals over spindles. Similarly, the backscattered ($\omega>\pi/2$) intensity can be given as the weighted integral over the apple $A_r$. We also consider the exterior problem, where we aim to reconstruct a function supported on the exterior of the unit ball from its weighted integrals over apples.

In section \ref{spindlesec} we introduce a new spindle transform for the monochromatic forward scatter problem and introduce a generalization of the spindle transform which gives the integrals of a function over the surfaces of revolution of a particular class of symmetric curves. We prove the injectivity of the generalized spindle transform on the domain of smooth functions compactly supported on the intersection of a hollow ball with the upper half space $x_3>0$. We show that our problem can be decomposed as a set of one dimensional inverse problems, which we then solve to provide an explicit expression for the harmonic coefficients of $f$.  In section \ref{TI} we introduce a new toric interior transform for the polychromatic forward scatter problem and prove its injectivity on the domain of smooth functions compactly supported on the intersection of a hollow ball and $x_3>0$. A new apple and apple interior transform are also introduced in section \ref{applesection} for the monochromatic and polychromatic backscatter problem. Their injectivity is proven on the domain of smooth functions compactly supported on the intersection of the exterior of the unit ball and $x_3>0$.

In section \ref{physics} we discuss possible approaches to the physical modelling of our problem for the case of a monochromatic and polychromatic photon source, and explain how this relates to the theory presented in section \ref{spindlesec}. 

In section \ref{results} we provide simulated density reconstructions of a test phantom via a discrete approach. We simulate data sets using the equations given in section \ref{spindlesec} and apply our reconstruction method with varying levels of added pseudo random noise. We also simulate the added effects due to the attenuation of the incoming and scattered rays in our data and see how this systematic error effects the quality of our reconstruction.

\section{A spindle transform}
\label{spindlesec}
We will now parameterize the set of points on a spindle $S_r$ in terms of spherical coordinates $(\rho,\theta,\varphi)$ and give some preliminary definitions before going on to define our spindle transform later in this section. 

\begin{figure}[!h]
\centering
\begin{tikzpicture}[scale=6]
\path ({-cos(30)},0.5) coordinate (S);
\path ({cos(30)}, 0.5) coordinate (D);
\path (-0.5,0.86) coordinate (w);
\path (-0.13,1.22) coordinate (a);
\draw [domain=30:150] plot ({cos(\x)}, {sin(\x)});
\draw [thin, dashed] (S) -- (D);
\draw [very thin] (S) -- (w);
\draw [very thin] (w) -- (D);
\draw [thin,dashed] (w)--(a);
\node at (-0.92,0.47) {$d$};
\node at (0.93,0.47) {$s$};
\node at (-0.41,0.88) {$\omega$};
\draw [ dashed] (0,0.5) circle [radius=0.87];
\draw [red] (-0.7,0.4) rectangle (-0.4,0.6);
\draw [->] (0,0.5)--(0,1.7);
\draw [->] (0,0.5)--(-1.2,0.5);
\node at (0.05,1.7) {$x_1$};
\node at (-1.2,0.45) {$x_3$};
\node at (-0.1,0.53) {$\varphi$};
\node at (0.4,1) {$C$};
\draw [->] (-0.9, 0.2)->(-0.7,0.4);
\node at (-0.9, 0.15) {object};
\draw [->] (-0.65,-0.45)->(-0.435,-0.253);
\node at (-0.65,-0.5) {scanning region};
\draw [<->] (0,0.5)--(0,0);
\node at (0.04,0.25) {$r$};
\draw [<->] (0,0)--(w);
\node at (-0.28,0.35) {$R$};
\draw [thin] (0,0.5)--(w);
\node at (-0.2,0.7) {$\rho$};
\end{tikzpicture}
\caption{A circular arc $C$ connecting source and detector points $s$ and $d$. The circular scanning region highlighted has unit radius.}
\label{fig3}
\end{figure}
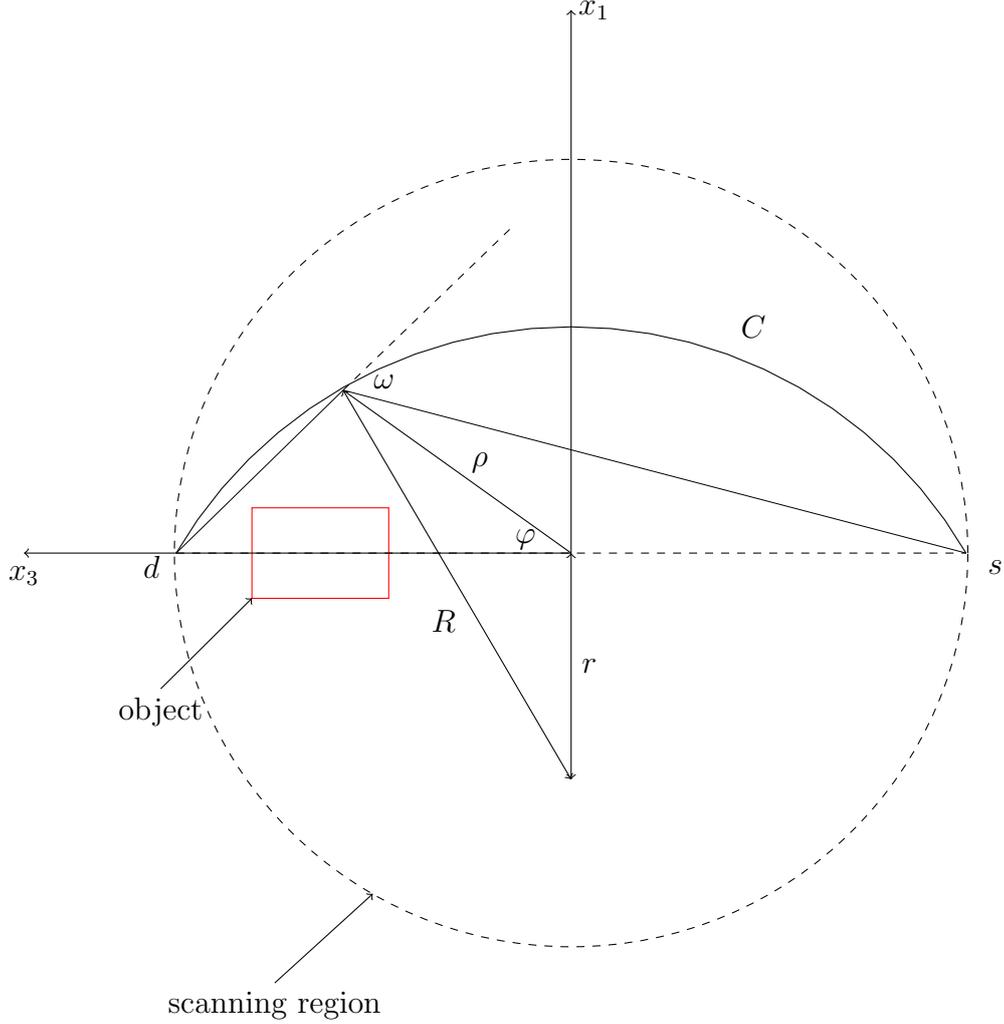

Consider the circular arc $C$ as illustrated in figure \ref{fig3}.
By the cosine rule, we have:
\begin{equation}
\label{cosine}
R^2=\rho^2+r^2+2r\rho\sin \varphi
\end{equation}
and hence:
\begin{equation}
\rho=\sqrt{r^2\sin^2\varphi+1}-r\sin\varphi.
\end{equation}
Let $B_{\epsilon_1,\epsilon_2}=\{x\in \mathbb{R}^3 : \epsilon_1 <|x|<\epsilon_2\}$ denote the set of points on a hollow ball with inner radius $\epsilon_1$ and outer radius $\epsilon_2$, and let $S^2$ denote the unit sphere. Let $Z=\mathbb{R}^{+}\times S^2$. For a function $f:\mathbb{R}^3\to \mathbb{R}$, let $F:Z\to \mathbb{R}$ be iits polar form $F(\rho,\theta,\varphi)=f(\rho\sin\varphi\cos\theta,\rho\sin\varphi\sin\theta,\rho\cos\varphi)$. We parameterize $h\in \text{SO}(3)$ in terms of Euler angles $\alpha$ and $\beta$, $h=U(\alpha)V(\beta)$, where $U$ and $V$ are rotations about the $x_3$ and $x_2$ axis respectively. We define an action of the rotation group $\text{SO}(3)$ on the set of real valued functions $f$ on $\mathbb{R}^3$ in the natural way, and define an action on the polar form as $(h\cdot F)=(h\cdot f)$. 

The arc element for the circular arc $C$ is given in \cite{NT} as:
\begin{equation}
\mathrm{d}v=\rho\sqrt{\frac{1+r^2}{1+r^2\sin^2{\varphi}}}\mathrm{d}\varphi
\end{equation}
and the area element for the spindle $S_r$ is:
\begin{equation}
\mathrm{d}A=\rho\sin\varphi\ \mathrm{d}v\mathrm{d}\theta.
\end{equation}
We define the spindle transform $\mathcal{S}:C_0^{\infty}(\mathbb{R}^3)\to C_0^{\infty}(Z)$ as:
\begin{equation}
\label{spindle}
\mathcal{S}f(r,\alpha,\beta)=\int_{0}^{2\pi}\int_{0}^{\pi}\rho^2\sin\varphi\sqrt{\frac{1+r^2}{1+r^2\sin^2{\varphi}}} (h\cdot F)\left(\rho,\theta,\varphi\right)\mid_{\rho=\sqrt{r^2\sin^2 \varphi +1}-r\sin \varphi}\mathrm{d}\varphi\mathrm{d}\theta
\end{equation}
where $h=U(\alpha)V(\beta)$.

This transform belongs to a larger class of integral transforms as we will now show. Let $p\in C^{1}([0,\frac{\pi}{2}])$ be a curve parameterized by a colatitude $\varphi \in [0,\frac{\pi}{2}]$. Then we define the class of symmetric curves $\rho \in C^{1}([0,1] \times [0,\pi])$ by:
\begin{equation}
\rho(r,\varphi)=p(\sin^{-1}(r\sin\varphi)).
\end{equation}
From this we define the generalized spindle transform $\mathcal{S}_{w,p}:C_0^{\infty}(\mathbb{R}^3)\to C_0^{\infty}([0,1]\times S^2)$ as:
\begin{equation}
\label{genspindle}
\mathcal{S}_{w,p}f(r,\alpha,\beta)=\int_{0}^{2\pi}\int_{0}^{\pi}w(r,\varphi)\rho\sin\varphi\sqrt{\rho^2+\left(\frac{\mathrm{d}\rho}{\mathrm{d}\varphi}\right)^2} (h\cdot F)\left(\rho,\theta,\varphi\right)\mid_{\rho=\rho(r,\varphi)}\mathrm{d}\varphi\mathrm{d}\theta.
\end{equation}
where the weighting $w$ is dependant only on $r$ and $\varphi$.

We now aim to prove injectivity of the generalized spindle transform $\mathcal{S}_{w,p}$ on the set of smooth functions on a hollow ball. First we give some definitions and background theory on spherical harmonics.

For integers $l\geq 0, \ |m|\leq l$, we define the spherical harmonics $Y_l^m$ as:
\begin{equation}
Y_l ^m(\theta,\varphi)=(-1)^m\sqrt{\frac{(2l+1)(l-m)!}{4\pi(l+m)!}}P_l ^m (\cos \varphi) e^{\mathrm{i}m\theta}
\end{equation}
where,
\begin{equation}
P_l^m(x)=(-1)^m(1-x^2)^{m/2}\frac{\mathrm{d}^m}{\mathrm{d}x^m}P_l(x)
\end{equation}
and
\begin{equation}
P_l(x)=\frac{1}{2^l}\sum_{k=0}^{l}\binom{l}{k}^2(x-1)^{l-k}(x+1)^{k},
\end{equation}
are Legendre polynomials of degree $l$. The spherical harmonics $Y_l^m$ form an orthonormal basis for $L^2(S^2)$, with the inner product:
\begin{equation}
\langle f,g\rangle=\int_{S^2}f\bar{g} \mathrm{d}\tau=\int_{0}^{\pi}\int_{0}^{2\pi}f(\theta,\varphi) \bar{g}(\theta,\varphi)\sin \varphi\ \mathrm{d}\theta \mathrm{d}\varphi.
\end{equation}
So we have:
\begin{equation}
\int_{S^2}Y_l ^m \bar{Y}_{l'} ^{m'} \mathrm{d}\tau=\int_{0}^{\pi}\int_{0}^{2\pi}Y_l ^m(\theta,\varphi) \bar{Y}_{l'} ^{m'}(\theta,\varphi)\sin \varphi\ \mathrm{d}\theta \mathrm{d}\varphi=\delta_{ll'}\delta_{mm'},
\end{equation}
where $\delta$ denotes the Kroneker delta.

From \cite{seeley} we have the theorem:
\begin{theorem}
\label{series}
Let $f\in C^{\infty}(\mathbb{R}^3)$ and let:
\begin{equation}
F_{lm}=\int_{S^2}F\bar{Y}_l ^m \mathrm{d}\tau,
\end{equation}
where $\mathrm{d}\tau$ is the surface measure on the sphere. Then the series:
\begin{equation}
F_N=\sum_{0\leq l\leq N} \sum_{|m|\leq l}F_{lm}Y_{l}^m
\end{equation}
converges uniformly absolutely on compact subsets of $Z$ to $F$.
\end{theorem}

We now show how the problem of reconstructing a density $f$ from its integrals $\mathcal{S}_{w,p} f$ for some curve $p\in C^{1}([0,\frac{\pi}{2}])$ can be broken down into a set of one dimensional inverse problems to solve for the harmonic coefficients $f_{lm}$:
\begin{lemma}
\label{mainlemma}
Let $f\in C_0^{\infty}(\mathbb{R}^3)$ and let $p\in C^{1}([0,\frac{\pi}{2}])$ be a curve, then:
\begin{equation}
\mathcal{S}_{w,p}f_{lm}(r)=2\pi\int_{0}^{\pi}w(r,\varphi)\rho\sin\varphi\sqrt{\rho^2+\left(\frac{\mathrm{d}\rho}{\mathrm{d}\varphi}\right)^2}F_{lm}(\rho)P_l(\cos\varphi)\mid_{\rho=\rho(r,\varphi)}\mathrm{d}\varphi
\end{equation}
where,
\begin{equation}
\mathcal{S}_{w,p}f_{lm}=\int_{S^2}(\mathcal{S}_{w,p}f)\bar{Y}_l ^m \mathrm{d}\tau
\end{equation}
and $P_l$ is a Legendre polynomial degree $l$.
\begin{proof}
Let $\tau\in S^2$ and let $h\in \text{SO}(3)$ act on a spherical harmonic $Y_l^m$ by $(h\cdot Y_l^m)(\tau)=Y_l^m(h\tau)$. Then we have:
\begin{equation}
\label{equ3}
\begin{split}
\mathcal{S}_{w,p}f(r,\alpha,\beta)&=\int_{0}^{2\pi}\int_{0}^{\pi}w(r,\varphi)\rho\sin\varphi\sqrt{\rho^2+\left(\frac{\mathrm{d}\rho}{\mathrm{d}\varphi}\right)^2} (h\cdot F)\left(\rho,\theta,\varphi\right)\mid_{\rho=\rho(r,\varphi)}\mathrm{d}\varphi\mathrm{d}\theta \\
&=\sum_{l\in \mathbb{N}}\sum_{|m|\leq l}\int_{0}^{2\pi}\int_{0}^{\pi}w(r,\varphi)\rho\sin\varphi\sqrt{\rho^2+\left(\frac{\mathrm{d}\rho}{\mathrm{d}\varphi}\right)^2} F_{lm}(\rho)(h\cdot Y_l^m)(\theta,\varphi)\mathrm{d}\varphi\mathrm{d}\theta
\end{split}
\end{equation}
where $h=U(\alpha)V(\beta)$.

We may write $h\cdot Y_l ^m$ as a linear combination of spherical harmonics of the same degree \cite[page 209]{sphharm}:
\begin{equation}
\label{equ27}
(h\cdot Y_l ^m)(\tau)=\sum_{|n|<l}Y_l ^n (\tau) D_{n,m} ^{(l)} (h^{-1})
\end{equation}
where the block diagonal entries $D_{n,m} ^{(l)}$ are defined:
\begin{equation}
D_{n,m} ^{(l)}(h)=D_{n,m} ^{(l)}(U(\alpha)V(\beta)U(\gamma))=e^{-i n\gamma}d_{n,m} ^{(l)}(\cos \beta) e^{-i m \alpha}
\end{equation}
and $d_{n,m} ^{(l)}$ is given, for $m=0$, by:
\begin{equation}
d_{n,0}^{(l)}(\cos \beta)=(-1)^n\sqrt{\frac{(l-n)!}{4\pi (l+n)!}}P_l ^n(\cos \beta).
\end{equation}
Here the $D^{l}$ are the irreducible blocks which form the regular representation of SO(3).

After the expansion (\ref{equ27}) is substituted into equation (\ref{equ3}), we see that the inserted sum is zero unless $n=0$ and we have:
\begin{equation}
\begin{split}
\label{equ8}
\mathcal{S}_{w,p}f(r,\alpha,\beta)&=2\pi \sum_{l\in \mathbb{N}} \sum_{|m|\leq l}\sqrt{\frac{(2l+1)}{4\pi}}D_{0,m} ^{(l)} (h^{-1})\int_{0}^{\pi}w(r,\varphi)\rho\sin\varphi\sqrt{\rho^2+\left(\frac{\mathrm{d}\rho}{\mathrm{d}\varphi}\right)^2} F_{lm}(\rho)P_l(\cos\varphi) \mathrm{d}\varphi.
\end{split}
\end{equation}
Since the regular representation of $\text{SO}(3)$ is unitary, we have $D_{0,m} ^{(l)} (h^{-1})=\bar{D}_{m,0} ^{(l)} (h)$ and from the above we have:
\begin{equation}
Y_l ^m (\alpha, \beta)=\sqrt{\frac{(2l+1)}{4\pi}}\bar{D}_{m,0} ^{(l)} (h).
\end{equation}
It follows that:
\begin{equation}
\begin{split}
\label{equ9}
\mathcal{S}_{w,p}f(r,\alpha,\beta)&=2\pi \sum_{l\in \mathbb{N}} \sum_{|m|\leq l}Y_l ^m (\alpha, \beta)\int_{0}^{\pi}w(r,\varphi) \rho\sin\varphi\sqrt{\rho^2+\left(\frac{\mathrm{d}\rho}{\mathrm{d}\varphi}\right)^2}F_{lm}(\rho)P_l(\cos\varphi) \mathrm{d}\varphi.
\end{split}
\end{equation}
From which we have:
\begin{equation}
\begin{split}
\mathcal{S}_{w,p}f_{lm}(r)&=\int_{S^2}\mathcal{S}_{w,p}f\bar{Y}_l ^m\mathrm{d}\tau\\
&=2\pi\int_{0}^{\pi}w(r,\varphi)\rho\sin\varphi\sqrt{\rho^2+\left(\frac{\mathrm{d}\rho}{\mathrm{d}\varphi}\right)^2}F_{lm}(\rho)P_l(\cos\varphi)\mid_{\rho=\rho(r,\varphi)}\mathrm{d}\varphi
\end{split}
\end{equation}
and the result is proven.
\end{proof}
\end{lemma}

We now go on to show that the set of one dimensional integral equations derived above are uniquely solvable for $f_{lm}$ for every $l\in \mathbb{N}$, $|m|\leq l$. Given the antisymmetry of the Legendre polynomials for odd $l$, we see that $\mathcal{S}_{w,p}f_{lm}=0$ for odd $l$, for any given $f\in C_0^{\infty}(\mathbb{R}^3)$ and curve $p\in C^{1}([0,\frac{\pi}{2}])$. So for the moment we focus on the reconstruction of the coefficients $f_{lm}$ for even $l$. 

We now have our main theorem:
\begin{theorem}
\label{maintheorem}
Let $p\in C^{1}([0,\frac{\pi}{2}])$ be such that $p(\varphi)\neq 0$ for all $\varphi \in [0,\frac{\pi}{2}]$, and let $g\in C^{1}([0,1))$ be defined as $g(x)=p(\sin^{-1}x)$. Let $A=\{\rho \in [0,1) : g(\rho)\in [0,\epsilon_1]\cup [\epsilon_2,\epsilon_3]\}$ and let $w : [0,1]\times [0,\pi] \to \mathbb{R}$ be a weighting such that $W(r,\rho)=w\left(r,\sin^{-1}(\rho/r)\right)$ and its first order partial derivative with respect to $r$ is continuous on $\{(r,\rho) : r\in A, \textnormal{min}(A)\leq\rho\leq r\}$ and $W(\rho,\rho)\neq 0$ on $A$. Let $f\in C_0^{\infty}(B_{0,\epsilon_1}\cup B_{\epsilon_2,\epsilon_3})$, where $\epsilon_1<p(0)<\epsilon_2$. Then $\mathcal{S}_{w,p}f$ determines the harmonic coefficients $F_{lm}$ uniquely for $\rho \in g([0,1))$ for all $l\in \{2k : k\in \mathbb{N}\}, |m|\leq l$.
\begin{proof}
Given the symmetry of the Legendre polynomials $P_l$ for even $l$, we have:
\begin{equation}
\label{equ4}
\begin{split}
\frac{1}{4\pi}\mathcal{S}_{w,p}f_{lm}(r)&=\int_{0}^{\frac{\pi}{2}}w(r,\varphi)\rho\sin\varphi\sqrt{\rho^2+\left(\frac{\mathrm{d}\rho}{\mathrm{d}\varphi}\right)^2}F_{lm}(\rho)P_l(\cos\varphi)\mid_{\rho=\rho(r,\varphi)}\mathrm{d}\varphi \\
&=\int_{0}^{\frac{\pi}{2}}w(r,\varphi)g(\rho)\sin\varphi\sqrt{g(\rho)^2+r^2\cos^2\varphi g'(\rho)^2}F_{lm}(g(\rho))P_l(\cos\varphi)\mid_{\rho=r\sin\varphi}\mathrm{d}\varphi.
\end{split}
\end{equation}
After making the substitution $
\rho=r\sin\varphi
$
, equation (\ref{equ4}) becomes:
\begin{equation}
\label{equ4.1}
\frac{1}{4\pi}\mathcal{S}_{w,p}f_{lm}(r)=\int_{0}^{r}\frac{F_{lm}(g(\rho))K_l(r,\rho)}{\sqrt{r-\rho}}\mathrm{d}\rho
\end{equation}
a Volterra integral equation of the first kind with weakly singular kernel, where:
\begin{equation}
K_l(r,\rho)=W(r,\rho)\frac{\rho g(\rho)\sqrt{g(\rho)^2+(r^2-\rho^2)g'(\rho)^2}}{r\sqrt{r+\rho}}P_l\left(\sqrt{1-\frac{\rho^2}{r^2}}\right).
\end{equation}
As $F_{lm}\circ g$ is zero for $\rho$ close to $0$, given our prior assumptions regarding $W$ and given that $g\in C^{1}([0,1))$, we can see that $K_l$ and its first order derivative with respect to $r$ is continuous on the support of $F_{lm}\circ g$.

Multiplying both sides of equation (\ref{equ4.1}) by $1/\sqrt{z-r}$ and integrating with respect to $r$ over the interval $[0,z]$, yields:
\begin{equation}
\label{equ7}
\frac{1}{4\pi}\int_{0}^{z}\frac{\mathcal{S}_{w,p}f_{lm}(r)}{\sqrt{z-r}}\mathrm{d}r=\int_{0}^{z}F_{lm}(g(\rho))\left[\int_{\rho}^{z}\frac{K_l(r,\rho)}{\sqrt{z-r}\sqrt{r-\rho}}\mathrm{d}r\right]\mathrm{d}\rho
\end{equation}
after changing the integration order. Making the substitution $r=\rho+(z-\rho)t$, gives:
\begin{equation}
\begin{split}
Q_{l}(z,\rho)&=\int_{\rho}^{z}\frac{K_l(r,\rho) }{\sqrt{z-r}\sqrt{r-\rho}}\mathrm{d}r \\
&=\int_{0}^{1}\frac{K_l(\rho+(z-\rho)t,\rho) }{\sqrt{t}\sqrt{1-t}}\mathrm{d}t
\end{split}
\end{equation}
from which we have: 
\begin{equation}
Q_l(z,z)=\pi K_l(z,z)=W(z,z)\frac{\pi c_l g(z)^2}{\sqrt{2z}}
\end{equation}
where,
\begin{equation}
c_l=P_l(0)=\frac{(-1)^{(l/2)}}{2^l}\binom{l}{l/2}.
\end{equation}
So by our assumptions that $p$ is non zero and $W$ is non zero on the diagonal, $Q_l(z,z)\neq 0$ on the support of $F_{lm}\circ g$.

Differentiating both sides of equation (\ref{equ7}) with respect to $z$ and rearranging gives:
\begin{equation}
\label{equ6}
v_{lm}(z)=\int_{0}^{z} F_{lm}(g(\rho))H_l(z,\rho)\mathrm{d}\rho+F_{lm}(g(z))
\end{equation}
which is a Volterra type integral equation of the second kind, where:
\begin{equation}
v_{lm}(z)=\frac{1}{4\pi^2 K_l(z,z)}\frac{\mathrm{d}}{\mathrm{d}z}\int_{0}^{z}\frac{\mathcal{S}_{w,p}f_{lm}(r)}{\sqrt{z-r}}\mathrm{d}r
\end{equation}
and,
\begin{equation}
H_l(z,\rho)=\frac{1}{\pi K_l(z,z)}\frac{\mathrm{d}}{\mathrm{d}z}Q_l(z,\rho).
\end{equation}
Given the continuity of $H_l$ on the support of $F_{lm}\circ g$ and the continuity of $F_{lm}\circ g$ on $[0,1]$, the Neumann series associated with equation (\ref{equ6}) converges and we may write our solution:
\begin{equation}
F_{lm}(g(z))=\int_{0}^{z}R_l(z,\rho) v_{lm}(\rho) \mathrm{d}\rho+v_{lm}(z)
\end{equation}
where the resolvent kernel,
\begin{equation}
R_l(z,\rho)=\sum_{i=1}^{\infty} H_{l,i} (z,\rho)
\end{equation}
is defined by,
\begin{equation}
H_{l,1} (z,\rho)=H_l(z,\rho), \ \ \ H_{l,i} (z,\rho) = \int_{\rho}^{z} H_l(z, x) H_{l,i-1} (x,\rho) \mathrm{d}x \ \ \forall i\geq 2.
\end{equation}
Since the series converges, the solution is unique and we may reconstruct $F_{lm}$ explicitly for $\rho\in g([0,1))$.
\end{proof}
\end{theorem}

Here as we can only reconstruct the coefficients $F_{lm}$ for even $l$, for a general $f\in C_0^{\infty}(\mathbb{R}^3)$ and curve $p\in C^{1}([0,\frac{\pi}{2}])$, the spindle data $\mathcal{S}_{w,p}f$ is insufficient to recover $f$ uniquely for $|x|\in p([0,\frac{\pi}{2}))$. However if we consider those functions $f\in C_0^{\infty}(\mathbb{R}^3)$ whose support lies in the upper half space $x_3>0$ we find that uniqueness is possible, as the following theorem shows:

\begin{theorem}
\label{coeff}
Let $U=\{(x_1,x_2,x_3)\in \mathbb{R}^3 : x_3>0\}$ and let $f\in C_0^{\infty}(U)$. Then the even coefficients $F_{lm}$ for $l\in \{2k : k\in \mathbb{N}\}=\mathbb{N}_e, |m|\leq l$ determine $F$ uniquely and:
\begin{equation}
F(\rho,\theta,\varphi)=2\sum_{l\in \mathbb{N}_e} \sum_{|m|\leq l}F_{lm}(\rho)Y_{l}^m (\theta,\varphi)
\end{equation}
for $\rho\geq 0$, $0\leq \theta\leq 2\pi$ and $0\leq\varphi\leq \frac{\pi}{2}$.
\begin{proof}
We can write $F$ as the sum of its even and odd coefficients:
\begin{equation}
F(\rho,\theta,\varphi)=F_{o}(\rho,\theta,\varphi)+F_{e}(\rho,\theta,\varphi)=\sum_{\text{odd}\ l\in \mathbb{N}} \sum_{|m|\leq l}F_{lm}(\rho)Y_{l}^m (\theta,\varphi)+\sum_{\text{even}\ l\in \mathbb{N}} \sum_{|m|\leq l}F_{lm}(\rho)Y_{l}^m (\theta,\varphi).
\end{equation}
Then from our assumption, we have:
\begin{equation}
\begin{split}
-F_e(\rho,\theta,\varphi)=F_o(\rho,\theta,\varphi)&=\sum_{\text{odd}\ l\in \mathbb{N}} \sum_{|m|\leq l}F_{lm}(\rho)Y_{l}^m (\theta,\varphi)\\
&=\sum_{m\in\mathbb{Z}}\left[\sum_{\text{odd}\ l\geq |m|}c(l,m)F_{lm}(\rho)P_l ^m(\cos \varphi ) \right] e^{im\theta}\\
&=\sum_{m\in\mathbb{Z}}F_o ^m(\rho,\varphi) e^{im\theta}
\end{split}
\end{equation}
for all $\rho\geq 0$, $0\leq \theta\leq 2\pi$ and $\frac{\pi}{2}\leq\varphi\leq \pi$, where:
\begin{equation}
c(l,m)=(-1)^m\sqrt{\frac{(2l+1)(l-m)!}{4\pi(l+m)!}}.
\end{equation}

We have:
\begin{equation}
\begin{split}
-F_e^m(\rho,\varphi)=-\frac{1}{2\pi}\int_{0}^{2\pi}F_e(\rho,\theta,\varphi) e^{-im\theta}\mathrm{d}\theta =F_o ^m(\rho,\varphi) \ \ \text{for}\ \ \rho\geq 0,\ \frac{\pi}{2}\leq\varphi\leq \pi.
\end{split}
\end{equation}
The associated Legendre polynomials $P_l ^m$ are symmetric when $l+m$ is even and antisymmetric otherwise. 
It follows that:
\begin{equation}
\begin{split}
F_o(\rho,\theta,\varphi)&=\sum_{m\in\mathbb{Z}}(-1)^{m+1}F_o ^m(\rho,\pi-\varphi) e^{im\theta}\\
&=\sum_{m\in\mathbb{Z}}(-1)^{m}F_e ^m(\rho,\pi-\varphi) e^{im\theta}\\
&=\sum_{m\in\mathbb{Z}}F_e ^m(\rho,\varphi) e^{im\theta}
\end{split}
\end{equation}
for $\rho\geq 0$, $0\leq \theta\leq 2\pi$ and $0\leq\varphi\leq \frac{\pi}{2}$. The result follows.
\end{proof}
\end{theorem}

\begin{corollary}
\label{corr0.5}
Let $f\in C_0^{\infty}(B_{\epsilon_1,\epsilon_2}\cap U)$ for some $0<\epsilon_1<\epsilon_2<1$. Let $\epsilon_1 \leq \epsilon\leq \epsilon_2$ and let $\delta=\frac{1-\epsilon^2}{2\epsilon}$. Then $\mathcal{S}f$ known for $0\leq r\leq \delta$ and for all $(\alpha,\beta)\in S^2$ determines $f$ uniquely for $\epsilon\leq |x|\leq 1$.
\begin{proof}
Let $p(\varphi)=\sqrt{1+\delta^2\sin^2\varphi}-\delta\sin\varphi$ and let $w\equiv 1$. Then $\mathcal{S}f(\delta r,\alpha,\beta)=\mathcal{S}_{w,p}f(r,\alpha,\beta)$ for $r\in [0,1]$, $(\alpha,\beta)\in S^2$. The result follows from Theorems \ref{maintheorem} and \ref{coeff}.
\end{proof}
\end{corollary}

\subsection{A toric interior transform}
\label{TI}
In the previous section we considered the three dimensional Compton scatter tomography problem for a monochromatic source and energy sensitive detector pair.  The polychromatic source case is covered here and we consider the modifications needed in our model to describe a full spectrum of initial photon energies.

In an X-ray tube electrons are accelerated by a large voltage ($E_m$keV) towards a target material and photons are emitted. Due to conservation of energy, the emitted photons have energy no greater than $E_m$keV. So for a given measured energy $E_s$, where $\frac{E_m}{1+2E_m/E_0}<E_s<E_m$, the set of scatterers is the union of spindle tori corresponding to scattering angles $\omega$ in the range $0<\omega<\cos^{-1}\left(1-\frac{E_0\left(E_m-E_s\right)}{E_{s}E_m}\right)$ (corresponding to energies $E$ in the range $E_s<E<E_m$). That is, the set of scatterers is a torus interior:
\begin{equation}
I_{r}=\left\{(x_1,x_2,x_3)\in \mathbb{R}^3 : \left(r-\sqrt{x_1^2+x_2^2}\right)^2+x_3^2<1+r^2\right\}
\end{equation}
where $r$ is determined by the scattered energy measured $E_s$. See \cite{Me} for an explanation of the two dimensional case.

With this in mind we define the spindle interior transform $\mathcal{I} : C_0^{\infty}(\mathbb{R}^3)\to C_0^{\infty}(Z)$ as:
\begin{equation}
\label{interior}
\mathcal{I}f(r,\alpha,\beta)=\int_{0}^{2\pi}\int_{0}^{\pi}\int_{0}^{\sqrt{r^2\sin^2\varphi+1}-r\sin\varphi}\rho^2\sin\varphi (h\cdot F)(\rho,\theta,\varphi)\mathrm{d}\rho\mathrm{d}\varphi\mathrm{d}\theta.
\end{equation}
and we have the uniqueness theorem for a weighted spindle interior transform:
\begin{theorem}
\label{interiortheorem}
Let $f\in C_0^{\infty}(B_{\epsilon_1,\epsilon_2}\cap U)$ for some $0<\epsilon_1<\epsilon_2<1$ and and let $\delta_2=\frac{1-\epsilon_2^2}{2\epsilon_2}$ and $\delta_1=\frac{1-\epsilon_1^2}{2\epsilon_1}$. Let $\tilde{f}$ be defined as:
\begin{equation}
\frac{1}{\sqrt{\frac{(1-|x|^2)^2}{4|x|^4}+1}}\tilde{f}(x)=\left(\frac{1-|x|^2}{2|x|^2}-\frac{\frac{(1-|x|^2)^2}{4|x|^4}}{\sqrt{\frac{(1-|x|^2)^2}{4|x|^4}+1}}\right) f\left(x\right).
\end{equation}
Define the weighted interior transform:
\label{TIW}
\begin{equation}
\label{polyspindle}
\begin{split}
\mathcal{I}_{w}f(r',\alpha,\beta)&=\int_{0}^{r'}\int_{0}^{2\pi}\int_{0}^{\pi}w(r',t, \varphi)\frac{\rho^2\sin\varphi}{t\sqrt{\frac{\sin^2\varphi}{t^2}+1}} (h\cdot \tilde{F})\left(\rho,\theta,\varphi\right)\mid_{\rho=\sqrt{\frac{\sin^2 \varphi}{t^2} +1}-\frac{\sin \varphi}{t}}\mathrm{d}\varphi\mathrm{d}\theta\mathrm{d}t,
\end{split}
\end{equation}
where $h=U(\alpha)V(\beta)$ and $r'=1/r$. Let us suppose that the weighting $w$ satisfies the following:
\begin{enumerate}
\item $w$ can be decomposed as $w(r',t, \varphi)=w_1(r',t)w_2(t,\varphi)$.
\item $W_2(t,\rho)=w_2\left(1/t,\sin^{-1}(\rho/t)\right)$ and its first order partial derivative with respect to $t$ is continuous on the triangle $T=\{(t,\rho)\in \mathbb{R}^2 : \epsilon_1\leq t\leq \epsilon_2, \epsilon_1\leq \rho\leq t\}$ and $W_2(\rho,\rho)\neq 0$ on $[\epsilon_1,\epsilon_2]$.
\item $w_1(r',t)$ and its first order partial derivatives are bounded on $0\leq r'<\infty$ and $w_1(r',r')\neq 0$ for $0\leq r'<\infty$.
\end{enumerate}
Then $\mathcal{I}_{w}f$ known for all $r'\in[0,\infty)$ and for all $(\alpha,\beta)\in S^2$ determines $f$ uniquely for $0<|x|<1$.

\begin{proof}
We have:
\begin{equation}
\label{equ20}
\mathcal{I}_{w}f(r',\alpha,\beta)=\int_{0}^{r'}w_1(r',t)G(t,\alpha,\beta)\mathrm{d}t,
\end{equation}
where,
\begin{equation}
G(t,\alpha,\beta)=\int_{0}^{2\pi}\int_{0}^{\pi}w_2(t, \varphi)\frac{\rho^2\sin\varphi}{t\sqrt{\frac{\sin^2\varphi}{t^2}+1}} (h\cdot \tilde{F})\left(\rho,\theta,\varphi\right)\mid_{\rho=\sqrt{\frac{\sin^2 \varphi}{t^2} +1}-\frac{\sin \varphi}{t}}\mathrm{d}\varphi\mathrm{d}\theta.
\end{equation}
Differentiating both sides of equation (\ref{equ20}) with respect to $r'$ and rearranging yields:
\begin{equation}
\label{equ21}
g(r',\alpha,\beta)=\int_{0}^{r'}L(r',t)G(t,\alpha,\beta)\mathrm{d}t+G(r',\alpha,\beta)
\end{equation}
where,
\begin{equation}
g(r',\alpha,\beta)=-\frac{\frac{\mathrm{d}}{\mathrm{d}r'} \mathcal{I}_{w}f(r',\alpha,\beta)}{w_1(r',r')}
\end{equation}
and,
\begin{equation}
L(r',t)=\frac{\frac{\mathrm{d}}{\mathrm{d}r'}w_1(r',t)}{w_1(r',r')}.
\end{equation}
Given our prior assumptions regarding $w_1$, we can solve the Volterra equation of the second kind (\ref{equ21}) uniquely for $G_{lm}(t)$ for $0< t<\infty$. From which we have:
\begin{equation}
\begin{split}
G\left(\frac{1}{\delta_1t},\alpha,\beta\right)&=\int_{0}^{2\pi}\int_{0}^{\pi}w_2\left(\frac{1}{\delta_1t}, \varphi\right)\frac{\delta_1t\rho^2\sin\varphi}{\sqrt{\delta^2_1t^2\sin^2\varphi+1}} (h\cdot \tilde{F})\left(\rho,\theta,\varphi\right)\mid_{\rho=\sqrt{\delta^2_2t^2\sin^2 \varphi} +1-\delta_2t\sin \varphi}\mathrm{d}\varphi\mathrm{d}\theta\\
&=\frac{\delta_1t}{\sqrt{1+\delta^2_1t^2}}\mathcal{S}_{w_3,p}\tilde{f}(t,\alpha,\beta)
\end{split}
\end{equation}
for $t\in (0,1]$, where $w_3(t,\varphi)=w_2\left(\frac{1}{\delta_1t}, \varphi\right)$ and $p(\varphi)=\sqrt{1+\delta_1^2\sin^2\varphi}-\delta_1\sin\varphi$. By our assumptions regarding $w_2$, the weighting $w_3$ satisfies the conditions of Theorem \ref{maintheorem}, as does the curve $p$. The result follows from Theorem \ref{maintheorem}.
\end{proof}
\end{theorem}
\begin{corollary}
Let $f\in C_0^{\infty}(B_{\epsilon_1,\epsilon_2}\cap U)$ for some $0<\epsilon_1<\epsilon_2<1$. Then $\mathcal{I}f$ as defined in equation (\ref{interior}) known for all $r\in[0,\infty)$ and for all $(\alpha,\beta)\in S^2$ determines $f$ uniquely for $0<|x|<1$.
\begin{proof}
Let $\tilde{f}$ be defined as in Theorem \ref{interiortheorem}. Then, after making the substitution $\rho=\sqrt{\frac{\sin^2 \varphi}{t^2} +1}-\frac{\sin \varphi}{t}$ in equation (\ref{interior}), we have:
\begin{equation}
\begin{split}
\mathcal{I}f(r,\alpha,\beta)&=\int_{0}^{\frac{1}{r}}\int_{0}^{2\pi}\int_{0}^{\pi}\frac{\rho^2\sin\varphi}{t\sqrt{\frac{\sin^2\varphi}{t^2}+1}} (h\cdot \tilde{F})\left(\rho,\theta,\varphi\right)\mid_{\rho=\sqrt{\frac{\sin^2 \varphi}{t^2} +1}-\frac{\sin \varphi}{t}}\mathrm{d}\varphi\mathrm{d}\theta\mathrm{d}t\\
&=\mathcal{I}_{w}f(1/r,\alpha,\beta)
\end{split}
\end{equation}
for all $r\in(0,\infty)$ and for all $(\alpha,\beta)\in S^2$, where the weighting $w\equiv 1$. The result follows from Theorem \ref{interiortheorem}.
\end{proof}
\end{corollary}

The advantage of using a polychromatic source (e.g. an X-ray tube) over a monochromatic source, which would most commonly be some type of gamma ray source, is that they have a significantly higher output intensity and so the data acquisition is faster. They are also safer to handle and use and are already in use in many fields of imaging. The downside, when compared to using a monochromatic source, would be a decrease in efficiency of our reconstruction algorithm and the added differentiation step required in the inversion process. This makes the polychromatic problem more ill posed, and so small errors in our measurements would be more greatly amplified in the reconstruction. We should consider both source types for further testing to determine an optimal imaging technique.

\subsection{The exterior problem}
\label{applesection}
Here we consider the exterior problem, and show that a full set of apple integrals (which represent the backscattered intensity) are sufficient to reconstruct a compactly supported density on the exterior of the unit ball.
\subsubsection{An apple transform}
For backscattered photons (for scattering angles $\omega>\pi/2$) the surface of scatterers is an apple $A_r$. Refer back to figure \ref{fig2}. The spherical coordinates $(\rho,\theta,\varphi)$ of the points on $A_r$ can be parameterized as follows:
\begin{equation}
\rho=\sqrt{r^2\sin^2\varphi +1}+r\sin\varphi,\ \ \ 0\leq \theta\leq 2\pi,\ \ \ 0\leq \varphi\leq \pi.
\end{equation}

We define the apple transform $\mathcal{A} : C_0^{\infty}(\mathbb{R}^3)\to C_0^{\infty}(Z)$ as:
\begin{equation}
\label{apple}
\begin{split}
\mathcal{A}f(r,\alpha,\beta)=\int_{0}^{2\pi}\int_{0}^{\pi}\rho^2\sin\varphi\sqrt{\frac{1+r^2}{1+r^2\sin^2{\varphi}}} (h\cdot F)(\rho,\theta,\varphi)\mid_{\rho=\sqrt{r^2\sin^2 \varphi +1}+r\sin \varphi}\mathrm{d}\varphi\mathrm{d}\theta,
\end{split}
\end{equation}
where $h=U(\alpha)V(\beta)$. We have the uniqueness theorem for the apple transform for functions supported on the exterior of the unit ball:
\begin{theorem}
\label{appletheorem}
Let $f\in C_0^{\infty}(B_{\epsilon_1,\epsilon_2}\cap U)$ for some $1<\epsilon_1<\epsilon_2<\infty$. Let $\epsilon_1 \leq \epsilon\leq \epsilon_2$ and let $\delta=\frac{\epsilon^2-1}{2\epsilon}$. Then $\mathcal{A}f$ known for $0\leq r\leq \delta$ and for all $(\alpha,\beta)\in S^2$ determines $f$ uniquely for $1\leq |x|\leq \epsilon$.
\begin{proof}
Let $p(\varphi)=\sqrt{1+\delta^2\sin^2\varphi}+\delta\sin\varphi$ and let $w\equiv 1$. Then $\mathcal{A}f(\delta r,\alpha,\beta)=\mathcal{S}_{w,p}f(r,\alpha,\beta)$ for $r\in [0,1]$, $(\alpha,\beta)\in S^2$. The result follows from Theorems \ref{maintheorem} and \ref{coeff}.
\end{proof}
\end{theorem}
\subsubsection{An apple interior transform}
Similar to our discussion at the start of section \ref{TI}, if the source is polychromatic the set of backscatterers is an apple interior. Hence we define the apple interior transform $\mathcal{AI} : C_0^{\infty}(\mathbb{R}^3)\to C_0^{\infty}(Z)$:
\begin{equation}
\label{appleinterior}
\mathcal{AI}f(r,\alpha,\beta)=\int_{0}^{2\pi}\int_{0}^{\pi}\int_{1}^{\sqrt{r^2\sin^2\varphi+1}+r\sin\varphi}\rho^2\sin\varphi (h\cdot F)(\rho,\theta,\varphi)\mathrm{d}\rho\mathrm{d}\varphi\mathrm{d}\theta.
\end{equation}
and we have a theorem for its injectivity:
\begin{theorem}
Let $f\in C_0^{\infty}(B_{\epsilon_1,\epsilon_2}\cap U)$ for some $1<\epsilon_1<\epsilon_2<\infty$. Let $\epsilon_1 \leq \epsilon\leq \epsilon_2$ and let $\delta=\frac{\epsilon^2-1}{2\epsilon}$. Then $\mathcal{AI}f$ known for $0\leq r\leq \delta$ and for all $(\alpha,\beta)\in S^2$ determines $f$ uniquely for $1\leq |x|\leq \epsilon$.
\begin{proof}
Let $\tilde{f}$ be defined as:
\begin{equation}
\frac{1}{\sqrt{\frac{(|x|^2-1)^2}{4|x|^4}+1}}\tilde{f}(x)=\left(\frac{\frac{(|x|^2-1)^2}{4|x|^4}}{\sqrt{\frac{(|x|^2-1)^2}{4|x|^4}+1}}+\frac{|x|^2-1}{2|x|^2}\right) f\left(x\right).
\end{equation}
Then by Leibniz rule we have:
\begin{equation}
\begin{split}
r\frac{\mathrm{d}}{\mathrm{d}r}\mathcal{AI}f(r,\alpha,\beta)&=\int_{0}^{2\pi}\int_{0}^{\pi}\frac{\rho^2\sin\varphi}{\sqrt{1+r^2\sin^2\varphi}} (h\cdot \tilde{F})(\rho,\theta,\varphi)\mid_{\rho=\sqrt{r^2\sin^2 \varphi +1}+r\sin \varphi}\mathrm{d}\varphi\mathrm{d}\theta\\
&=\frac{1}{\sqrt{1+r^2}}\mathcal{A}\tilde{f}(r,\alpha,\beta)
\end{split}
\end{equation}
So $\mathcal{AI}f$ known for $0\leq r\leq \delta$ and for $(\alpha,\beta)\in S^2$ determines $\mathcal{A}\tilde{f}$ on the same set after differentiating. The result follows from Theorem \ref{appletheorem}.
\end{proof}
\end{theorem}

\section{A physical model}
\label{physics}
Here we explain how the theory presented in the previous section relates to what we measure in a practical setting. 

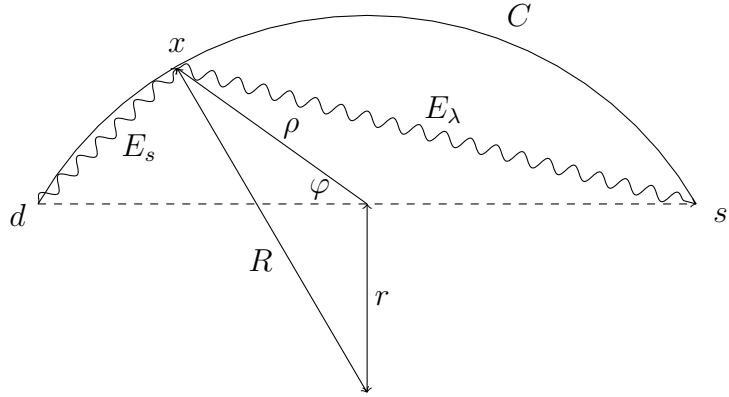
\begin{figure}[!h]
\centering
\begin{tikzpicture}[scale=5]
\path ({-cos(30)},0.5) coordinate (S);
\path ({cos(30)}, 0.5) coordinate (D);
\path (-0.5,0.86) coordinate (w);
\path (-0.13,1.22) coordinate (a);
\draw [domain=30:150] plot ({cos(\x)}, {sin(\x)});
\draw [thin, dashed] (S) -- (D);
\draw [snake it, ->] (S) -- (w);
\draw [snake it, ->] (w) -- (D);
\node at (-0.92,0.47) {$d$};
\node at (0.93,0.47) {$s$};
\node at (-0.125,0.54) {$\varphi$};
\node at (0.2,0.75) {$E_{\lambda}$};
\node at (-0.6,0.65) {$E_s$};
\node at (-0.5, 0.92) {$x$};
\node at (0.4,1) {$C$};
\draw [<->] (0,0.5)--(0,0);
\node at (0.04,0.25) {$r$};
\draw [<->] (0,0)--(w);
\node at (-0.28,0.35) {$R$};
\draw [thin] (0,0.5)--(w);
\node at (-0.2,0.7) {$\rho$};
\end{tikzpicture}
\caption{A scattering event occurs on the circular arc $C$ with initial photon energy $E_{\lambda}$ from the source $s$ to the detector $d$ with energy $E_s$.}
\label{fig5}
\end{figure}

We consider an intensity of photons scattering from a point $x$ as illustrated in figure \ref{fig5}.
The points $s$ and $d$ are the centre points of the source and detector respectively. The intensity of photons scattered from $x$ to $d$ with energy $E_s$ is:
\begin{equation}
\label{equ17}
\begin{split}
I\left(E_s\right)=I_0\left(E_{\lambda}\right)e^{-\int_{L_{sx}}\mu(E_{\lambda},Z)}
f\left(x\right)\mathrm{d}V\times \frac{\mathrm{d}\sigma}{\mathrm{d}\Omega}\left(E_s,\omega\right)e^{-\int_{L_{xd}}\mu(E_s,Z)}S\left(q,Z\right) 
\mathrm{d}\Omega.
\end{split}
\end{equation}
where $Z$ denotes the atomic number, $\omega$ the scattering angle and $L_{sx}$ and $L_{xd}$ are the line segments connecting $s$ to $x$ and $x$ to $d$ respectively. $f(x)$ denotes the electron density (number of electrons per unit volume) at the scattering point $x$ and $\mathrm{d}V$ is the volume measure. $f$ is the quantity to be reconstructed. 

Let $W_k(E_{\lambda})$ denote the incident photon flux (number of photons per unit area per unit time), energy $E_{\lambda}$, measured at a fixed distance $D$ from the source. Then the incident photon intensity $I_0$ can be written:
\begin{equation}
I_0(E_{\lambda})=\frac{tD^2W_k(E_{\lambda})}{(\rho\cos\varphi+1)^2+\rho^2\sin^2\varphi},
\end{equation}
where $t$ is the emission time.

The Klein-Nishina differential cross section $\mathrm{d}\sigma/\mathrm{d}\Omega$, is defined by:\begin{equation}\frac{\mathrm{d}\sigma}{\mathrm{d}\Omega}\left(E_s,\omega\right)=\frac{r_0^{2}}{2}{\left(\frac{E_s}{E_{\lambda}}\right)}^{2}\left(\frac{E_s}{E_{\lambda}}+\frac{E_{\lambda}}{E_s}-1+\cos^{2} \omega\right),\end{equation}
where $r_0$ is the classical electron radius and $\cos\omega=r/\sqrt{1+r^2}$. This predicts the scattering distribution for a photon off a free electron at rest. Given that the atomic electrons typically are neither free nor at rest, a correction factor is included, namely the incoherent scattering function $S\left(q,Z\right)$. Here $q=\frac{E_{\lambda}}{hc}\sin \left(\omega /2\right)$ is the momentum transferred by a photon with initial energy:
\begin{equation}
\label{inE}
E_{\lambda}=\frac{E_s}{1-\left(E_s/E_0\right)\left(1-\cos\omega\right)}
\end{equation}
scattering at an angle $\omega$, where $h$ is Planck's constant and $c$ is the speed of light. As the scattering function $S$ is dependant on the atomic number $Z$, we set $Z=Z_{\text{avg}}$ to some average atomic number as an approximation and interpolate values of $S(q,Z_{\text{avg}})$ from the tables given in \cite{hub}.

The solid angle subtended by $x$ and $d$ may be approximated as:
\begin{equation}\mathrm{d}\Omega=\frac{A}{4\pi}\times\frac{1-\rho\cos\varphi}{(1-\rho\cos\varphi)^2+\rho^2\sin^2\varphi},\end{equation}
where $A$ is the detector area. 

The exponential terms in equation (\ref{equ17}) account for the attenuation of the incoming and scattered rays. We approximate:
\begin{equation}
e^{-\int_{L_{sx}}\mu(E_{\lambda},Z)}e^{-\int_{L_{xd}}\mu(E_s,Z)}\approx 1
\end{equation}
In general this approximation is unrealistic. For example in medical CT, if we were scanning a relatively large (e.g. the size of someone's head) mass of organic material at a low energy (e.g. 50keV), the absorption would play a significant role and the above model would over approximate the data. However in an application where the objects are smaller (centimetres in diameter) and where we can scan at a higher energy ($\approx1$MeV), the effects due to absorption would be less prevalent and Compton scattering would be the dominant interaction. For example, in airport baggage screening typical hand luggage would be a small bag containing a few low effective $Z$ densities (e.g. nail varnish (Acetone), water, some plastic (polyethylene) etc.) and the rest may be clothes or air. Here also, as we are not worried about dosage, we can scan at high energies (e.g. using a high voltage X-ray tube or high emission energy gamma ray source). We will simulate the error due to attenuation later, in section \ref{results}, and show the effects of neglecting the attenuation in our reconstruction.

\subsubsection{The monochromatic case}
Let our density $f$ be supported on $B_{\epsilon_1,\epsilon_2}\cap U$ for some $0<\epsilon_1<\epsilon_2<1$ and let $\delta_2=\frac{1-\epsilon^2_2}{2\epsilon_2}$. If the source $s$ is monochromatic ($E_{\lambda}$ remains fixed), then the forward Compton scattered intensity measured is:
\begin{equation}
\label{equmodel}
I(r,\alpha,\beta)={ctD^2W_k(E_{\lambda})}\frac{\mathrm{d}\sigma_c}{\mathrm{d}\Omega}(r)\mathcal{S}_{w,p}f(r,\alpha,\beta),
\end{equation}
where $p\in C^1([0,\pi])$ is defined by $p(\varphi)=\sqrt{1+\delta^2_2\sin^2\varphi}-\delta_2\sin \varphi$, $c$ is a constant thickness and:
\begin{equation}
\frac{\mathrm{d}\sigma_c}{\mathrm{d}\Omega}(r)=\frac{\mathrm{d}\sigma}{\mathrm{d}\Omega}\left(E_s,\omega\right)S(q,Z_{\text{avg}})
\end{equation}
depends only on $r$. Here the variable $r\in [0,1]$ determines the  scattered energy $E_s$ and $(\alpha,\beta)\in S^2$ determines the source and detector position. The weighting $w$ is given by:
\begin{equation}
\label{weight}
w(r,\varphi)=\frac{A(1-\rho\cos\varphi)}{4\pi\left[(1-\rho\cos\varphi)^2+\rho^2\sin^2\varphi\right]\left[(\rho\cos\varphi+1)^2+\rho^2\sin^2\varphi\right]}.
\end{equation}
where $\rho=\sqrt{1+r^2\sin^2\varphi}-r\sin\varphi$. It is left to the reader to show that the weighting $w$ satisfies the conditions given in Theorem \ref{maintheorem}. After dividing through by the physical modelling terms in equation (\ref{equmodel}), we can invert the weighted spindle transform $\mathcal{S}_{w,p}f$ as in Theorem \ref{maintheorem} to obtain an analytic expression for $f$ in terms of the Compton scattered intensity $I$.
\subsubsection{The polychromatic case}
\begin{figure}[!h]
\centering
\includegraphics[scale=0.5]{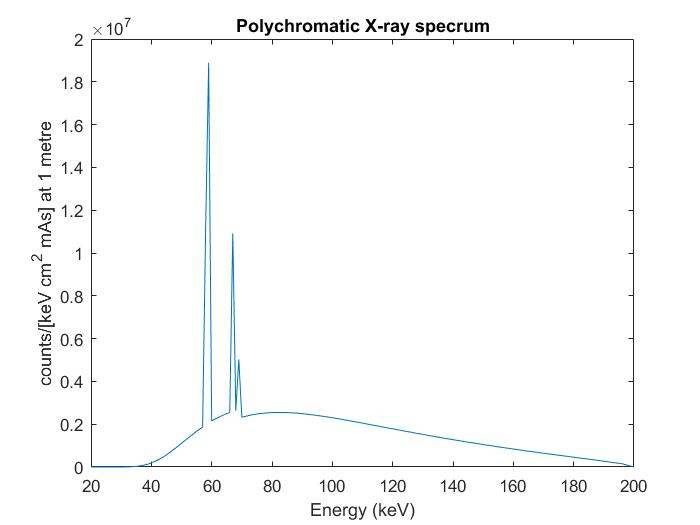}
\caption{A polychromatic Tungsten target spectrum with an X-ray tube voltage of $V=200$keV and a 1mm Copper filter. Spectrum calculated using SpekCalc \cite{spekcalc}.}
\label{F25}
\end{figure}

Here we have a spectrum of incoming photon energies. See figure \ref{F25}.
The Compton scattered intensity measured for a tube centre offset $r=1/r'$ for a source and detector position $(\alpha,\beta)\in S^2$ may be written:
\begin{equation}
I(r',\alpha,\beta)=stD^2\mathcal{I}_wf(r',\alpha,\beta),
\end{equation}
with the weighting:
\begin{equation}
w(r',t,\varphi)=\left[W_k(r',t)\frac{\mathrm{d}\sigma_c}{\mathrm{d}\Omega}(r',t)\right]w_2(t,\varphi)=w_1(r',t)w_2(t,\varphi),
\end{equation}
where $w_2=w$ as in equation (\ref{weight}). Here the scattering probabilty $\frac{\mathrm{d}\sigma_c}{\mathrm{d}\Omega}$ and the photon flux $W_k$ are dependant on $r'$ and the integration variable $t$ as in subsection \ref{TI}. While the weighting $w$ is separable as in Theorem \ref{TIW}, the incident photon flux $W_k(r',r')=W_k(E_m)=0$ for all $r'$, where $E_m$ is the maximum spectrum energy. Hence $w_1(r',r')\equiv 0$ and $w$ fails to meet the conditions of Theorem \ref{TIW}. To deal with this, let:
\begin{equation}
\label{approx}
w_{\text{avg}}(r')=\frac{1}{r'}\int_{0}^{r'}w_1(r',t)\mathrm{d}t
\end{equation}
and let $w_{\text{app}}(r',t,\varphi)=w_{\text{avg}}(r')w_2(t,\varphi)$. Then, although we sacrifice some accuracy in our forward model, $w_{\text{app}}$ would satisfy the conditions given in Theorem \ref{TIW} and we can obtain an analytic expression for the density. The error in our approximation is bounded by:
\begin{equation}
\left|w_{\text{avg}}(r')-w_1(r',t)\right|\leq \left|\max_{t\in [0,r']}w_1(r',t)-\min_{t\in [0,r']}w_1(r',t)\right|
\end{equation}
for all $r'\geq 0$, $0\leq t\leq r'$. So provided the changes in the incident photon energy $E_s$ ($E_s$ is determined by $t$) are negligible over the range $t\in [0,r']$, the error in $w_{\text{app}}$ would be small.
\section{Simulations}
\label{results}
Here we provide reconstructions of a test phantom density at a low resolution using simulated datasets of the unweighted spindle and spindle interior transforms, and simulate noise as additive psuedo Gaussian noise. To simulate data for each transform we discretize the integrals in equations (\ref{spindle}) and (\ref{interior}), and solve the least squares problem:
\begin{equation}
\label{tik}
\argmin_x\|Ax-b\|^2+\lambda^2\|x\|^2
\end{equation}
for some regularisation parameter $\lambda>0$, where $A$ is the discrete forward operator of the transform considered, $x$ is the vector of density pixel values and $b=Ax$ is the simulated transform data. We simulate perturbed data $b^{\epsilon}$ via:
\begin{equation}
\label{error}
b^{\epsilon}=b+\epsilon \frac{G\|b\|}{\sqrt{n}}
\end{equation}
where $G$ is a pseudo random vector of samples from the standard Gaussian distibution and $n$ is the number of entries in $b$. The proposed noise model has the property that $\|b-b^{\epsilon}\|/\|b\|\approx \epsilon$, so a noise level of $\epsilon$ is $\epsilon\times 100\%$ relative error.

\begin{figure}[h]
\hspace*{-2cm}
\begin{subfigure}{0.5\textwidth}
\includegraphics[width=1.2\linewidth, height=10.5cm]{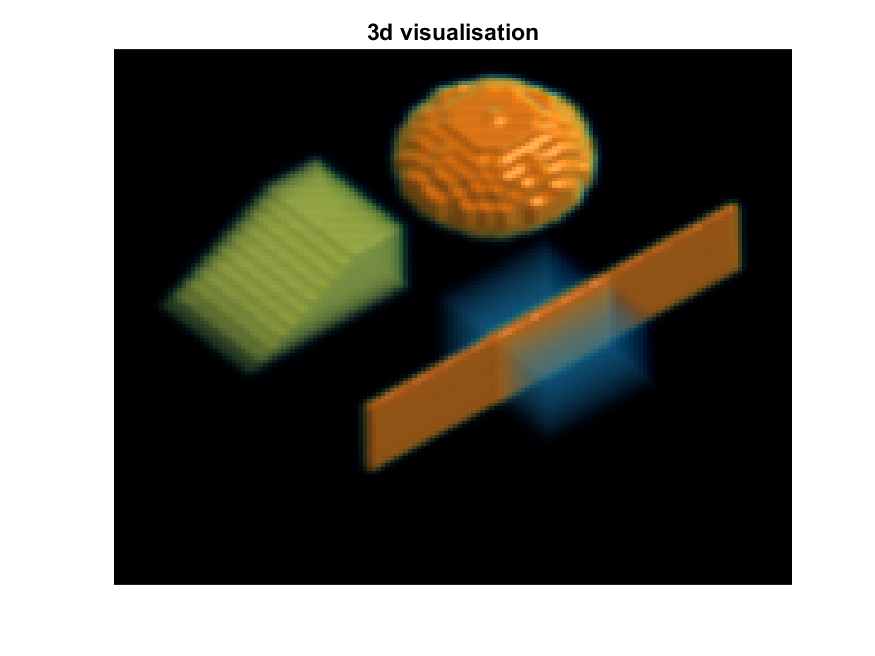} 
\end{subfigure}\hspace{15mm}
\begin{subfigure}{0.5\textwidth}
\includegraphics[width=1.2\linewidth, height=10.5cm]{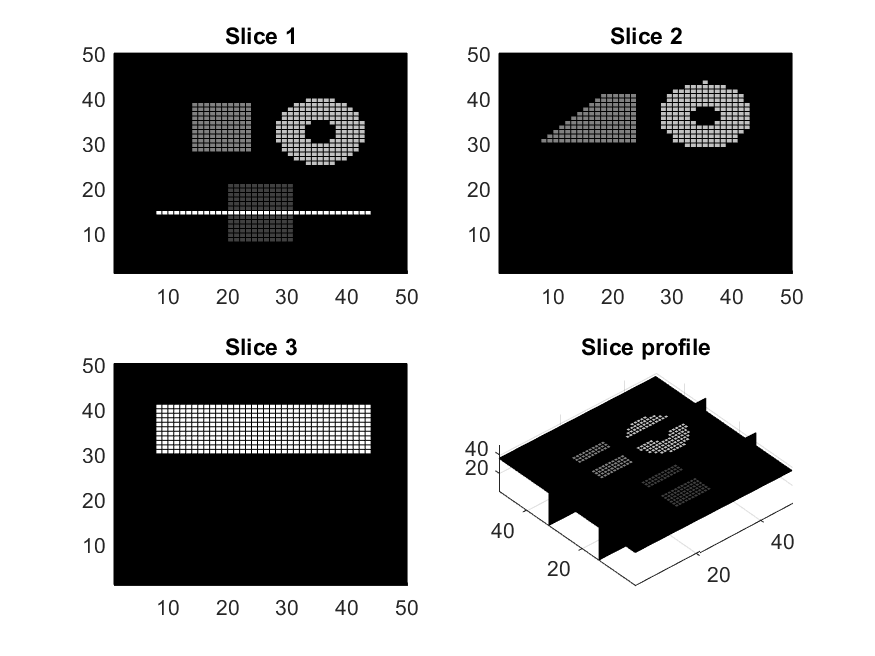}
\end{subfigure}
\caption{Test phantom.}
\label{F1}
\hspace*{-2cm}
\begin{subfigure}{0.5\textwidth}
\includegraphics[width=1.2\linewidth, height=10.5cm]{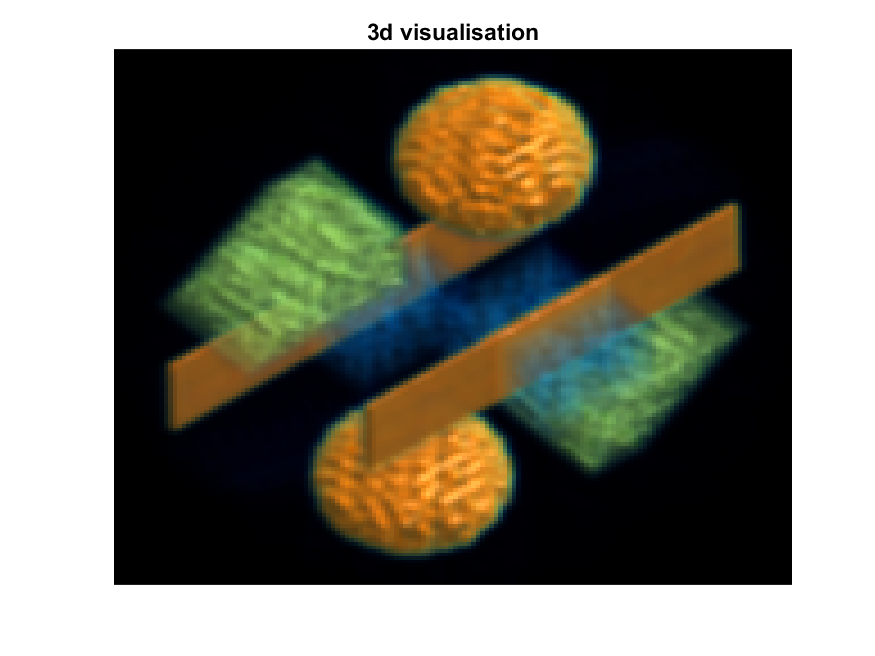} 
\end{subfigure}\hspace{15mm}
\begin{subfigure}{0.5\textwidth}
\includegraphics[width=1.2\linewidth, height=10.5cm]{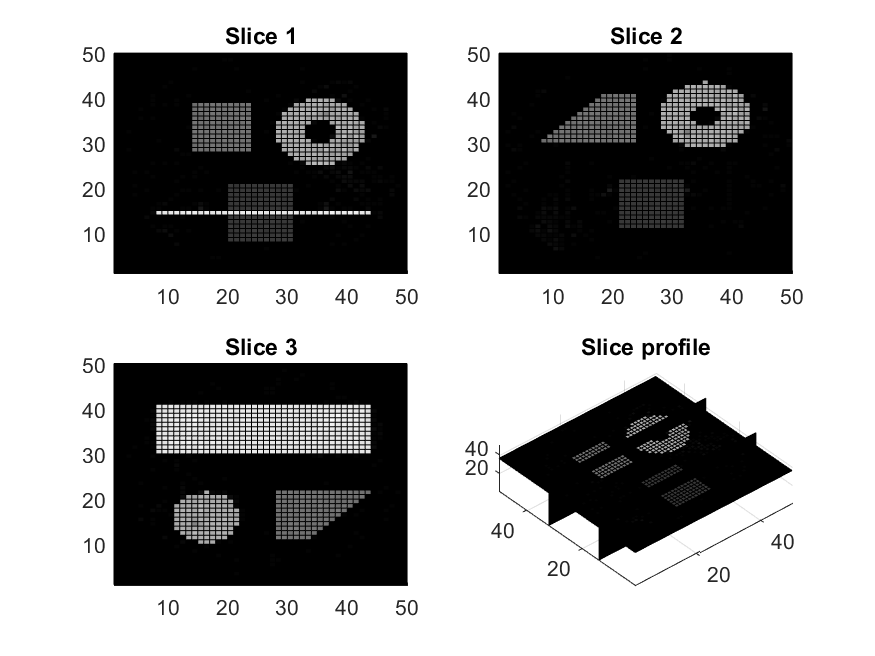}
\end{subfigure}
\caption{Test phantom reconstruction from spindle transform data no noise, 2000 iterations.}
\label{F2}
\end{figure}

We consider the test phantom displayed in figure \ref{F1}. The unit cube is discretised into $50\times50\times50$ pixels and a hollow ball, some stairs and a low density block with a metal sheet passing through it are placed in upper unit hemisphere. Slice profiles are given in the right hand figure to display the metal sheet and the hollow shell. We sample 25 $r$ values (corresponding to spindles with heights $H\in \{0.02+0.04i : 0\leq i <25\}$), 45 $\alpha$ values $\alpha\in \{\frac{2 \pi i}{45} : 0\leq i < 45\}$ and 45 $\beta$ values $\beta\in \{\frac{\pi}{90}(1+i) : 0\leq i<45\}$. So we have an underdetermined sparse system matrix $A$ with 50625 rows and 125000 columns. To reconstruct we solve the regularised problem (\ref{tik}) using the conjugate gradient least squares algorithm (CGLS) and pick our regularisation parameter $\lambda$ via a manual approach. In each reconstruction presented, any negative values are set to zero and the iteration number and noise level are given in the figure caption.

\begin{figure}[h]
\hspace*{-2cm}
\begin{subfigure}{0.5\textwidth}
\includegraphics[width=1.2\linewidth, height=10.5cm]{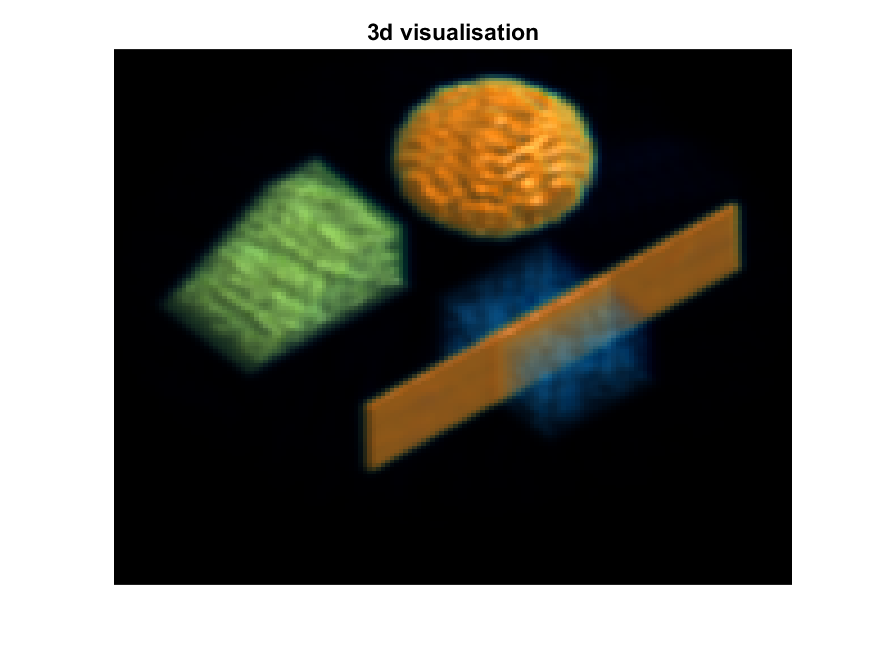} 
\end{subfigure}\hspace{15mm}
\begin{subfigure}{0.5\textwidth}
\includegraphics[width=1.2\linewidth, height=10.5cm]{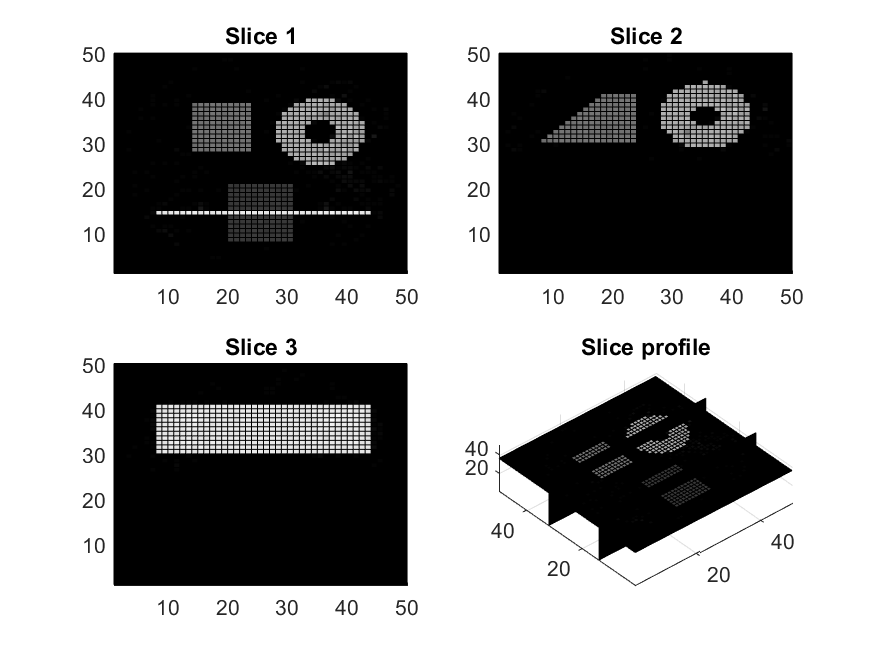}
\end{subfigure}
\caption{Test phantom reconstruction from spindle transform data no noise, set to zero in lower hemisphere, 2000 iterations.}
\label{F3}
\hspace*{-2cm}
\begin{subfigure}{0.5\textwidth}
\includegraphics[width=1.2\linewidth, height=10.5cm]{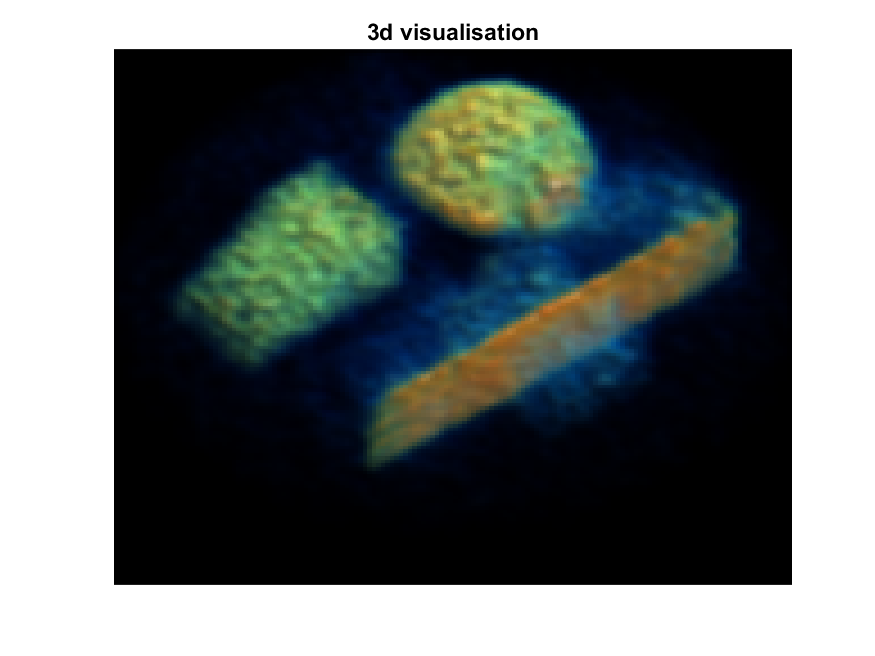} 
\end{subfigure}\hspace{15mm}
\begin{subfigure}{0.5\textwidth}
\includegraphics[width=1.2\linewidth, height=10.5cm]{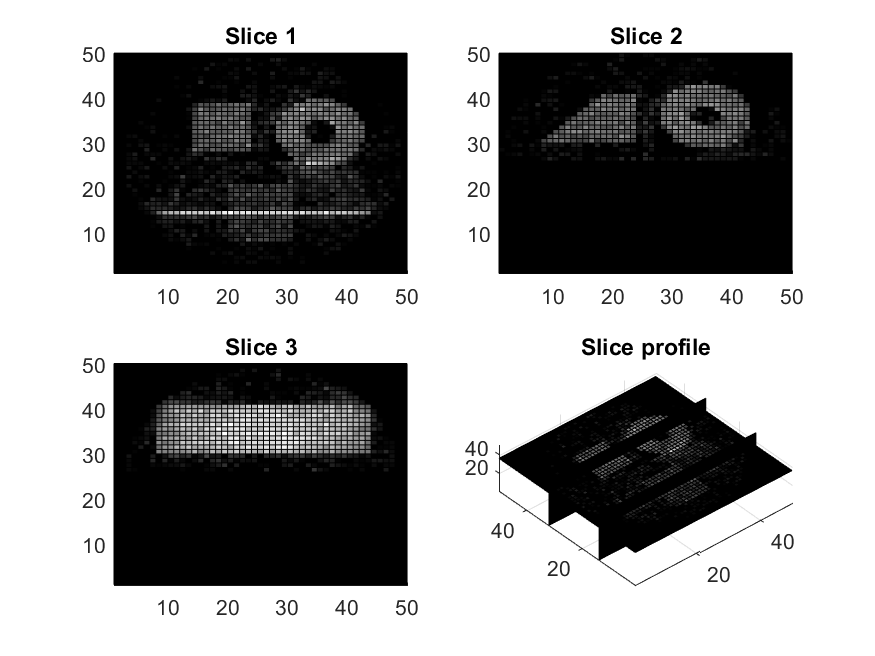}
\end{subfigure}
\caption{Test phantom reconstruction from spindle transform data, noise level $\epsilon=0.01$, 2000 iterations.}
\label{F4}
\end{figure}

In figure \ref{F2} we have presented a reconstruction of the test phantom in the absence of added noise. Given the symmetries involved in our geometry, the density has been rotated and reflected in the $xy$ plane in the reconstruction. To better visualise the reconstruction we set our reconstructed images to 0 in the lower half space. In figure \ref{F3} we display the same zero noise reconstruction but with the voxel values set to zero in the lower hemisphere. 

\begin{figure}[h]
\hspace*{-2cm}
\begin{subfigure}{0.5\textwidth}
\includegraphics[width=1.2\linewidth, height=10.5cm]{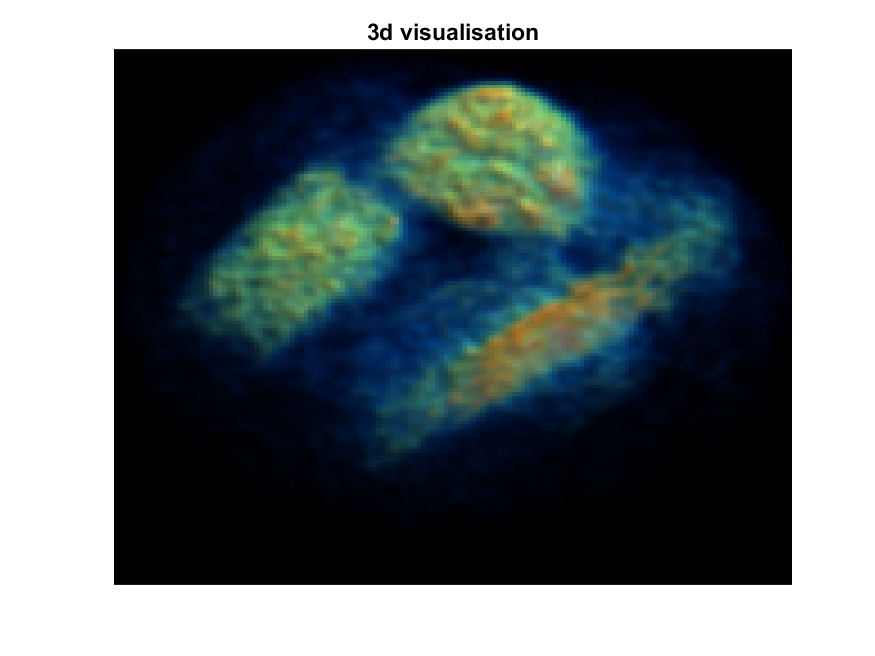} 
\end{subfigure}\hspace{15mm}
\begin{subfigure}{0.5\textwidth}
\includegraphics[width=1.2\linewidth, height=10.5cm]{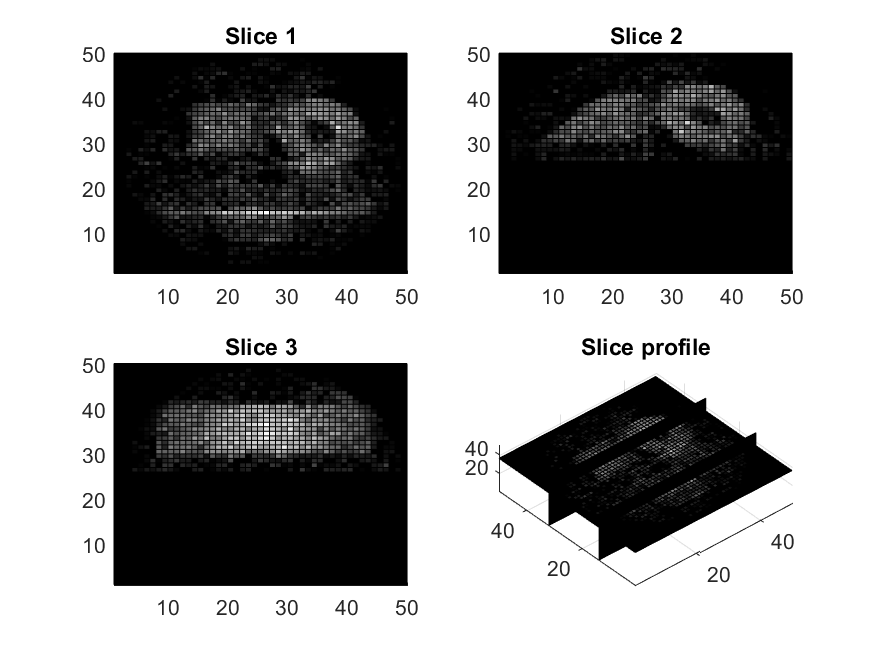}
\end{subfigure}
\caption{Test phantom reconstruction from spindle transform data, noise level $\epsilon=0.05$, 2000 iterations.}
\label{F5}
\hspace*{-2cm}
\begin{subfigure}{0.5\textwidth}
\includegraphics[width=1.2\linewidth, height=10.5cm]{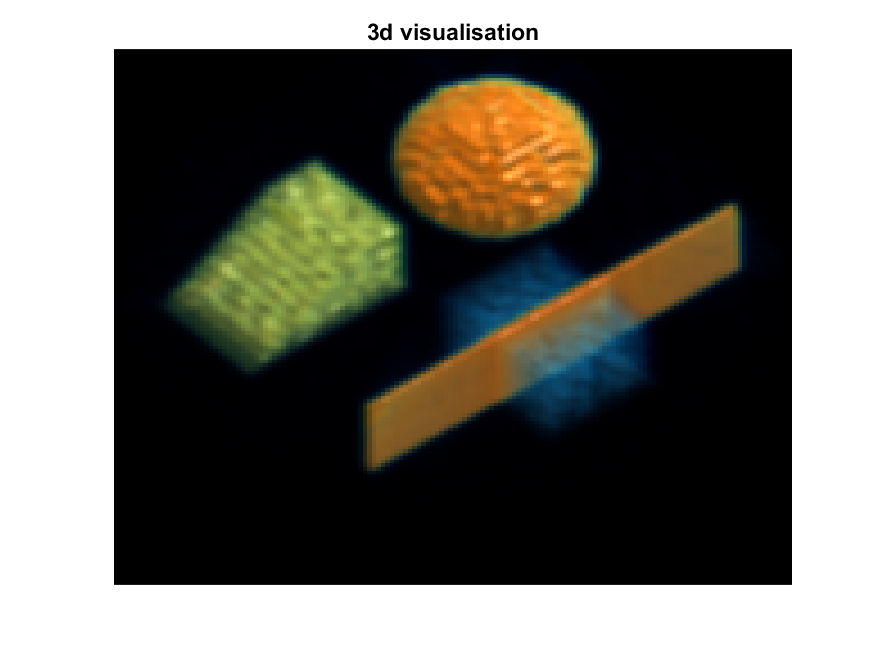} 
\end{subfigure}\hspace{15mm}
\begin{subfigure}{0.5\textwidth}
\includegraphics[width=1.2\linewidth, height=10.5cm]{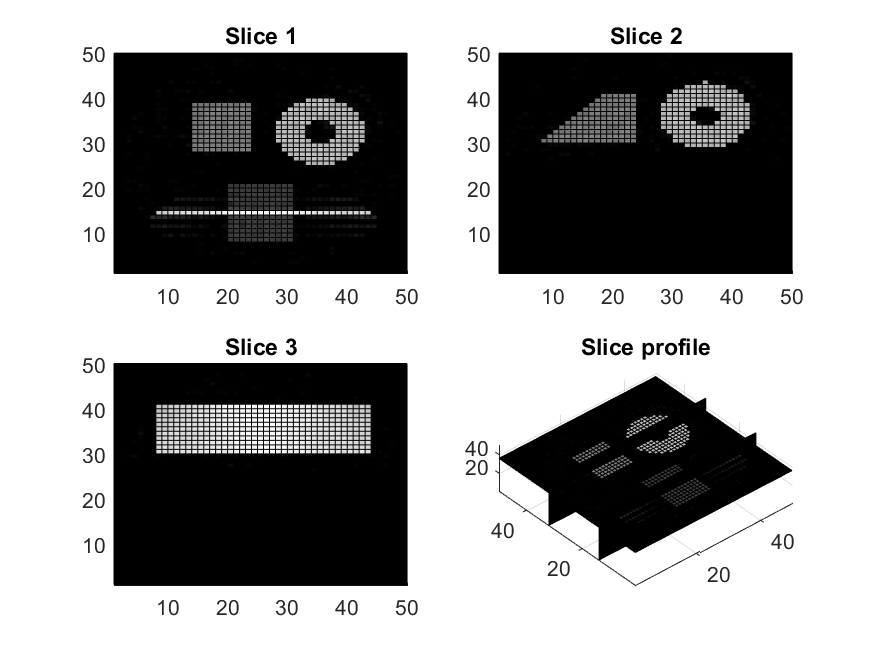}
\end{subfigure}
\caption{Test phantom reconstruction from spindle interior data no noise, 2000 iterations.}
\label{F6}
\end{figure}

\begin{figure}[h]
\hspace*{-2cm}
\begin{subfigure}{0.5\textwidth}
\includegraphics[width=1.2\linewidth, height=10.5cm]{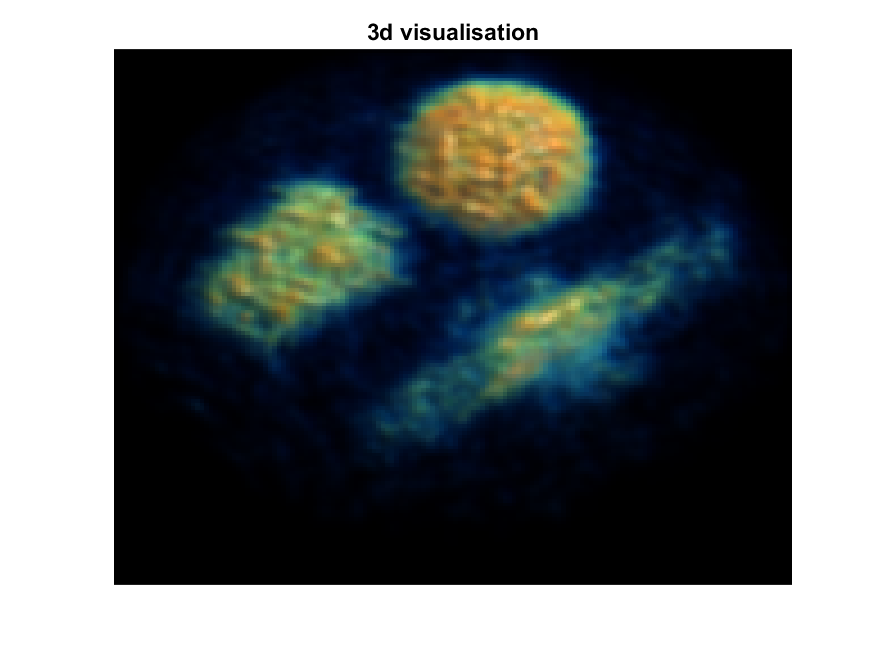} 
\end{subfigure}\hspace{15mm}
\begin{subfigure}{0.5\textwidth}
\includegraphics[width=1.2\linewidth, height=10.5cm]{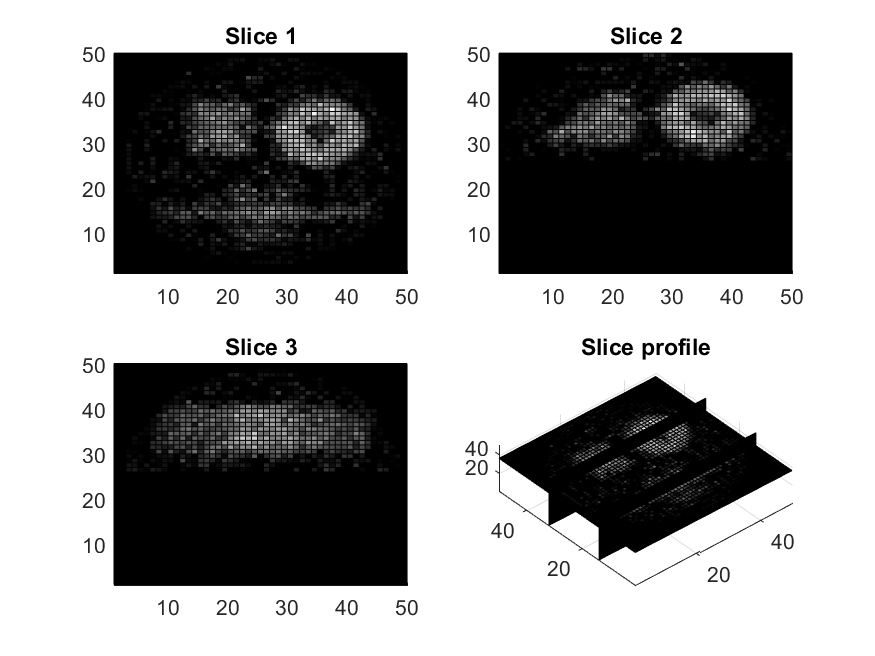}
\end{subfigure}
\caption{Test phantom reconstruction from spindle interior data, noise level $\epsilon=0.01$, 2000 iterations.}
\label{F7}
\hspace*{-2cm}
\begin{subfigure}{0.5\textwidth}
\includegraphics[width=1.2\linewidth, height=10.5cm]{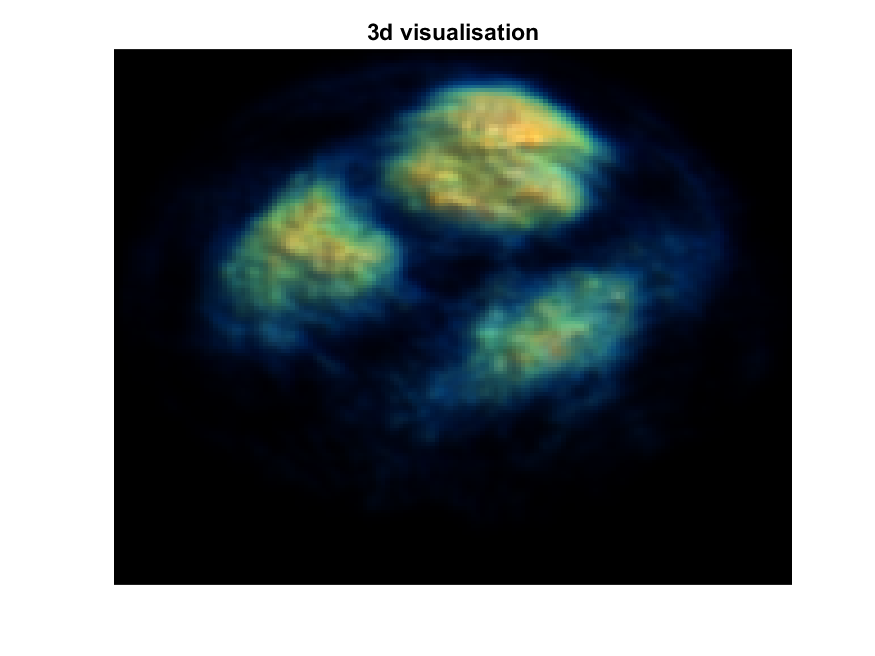} 
\end{subfigure}\hspace{15mm}
\begin{subfigure}{0.5\textwidth}
\includegraphics[width=1.2\linewidth, height=10.5cm]{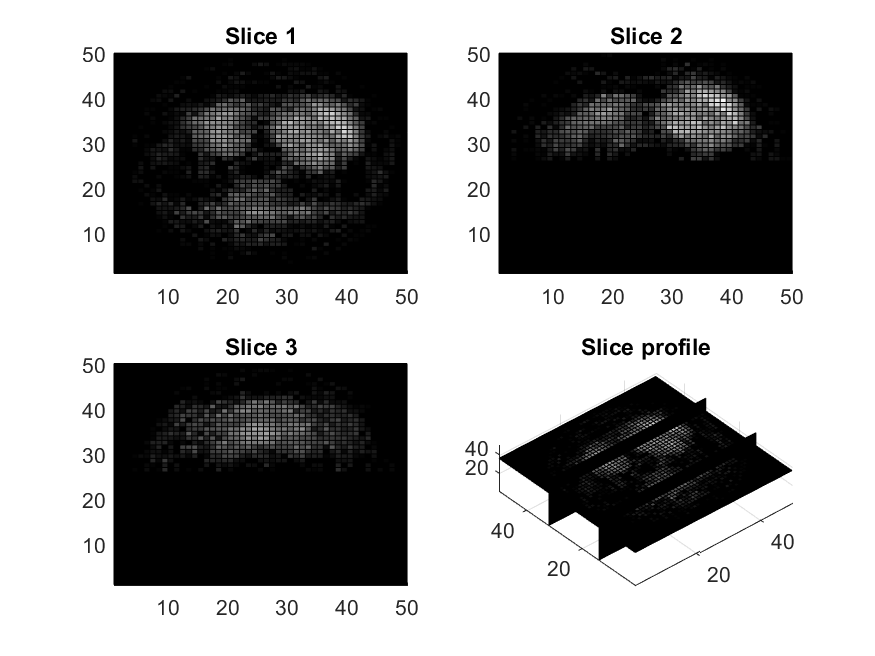}
\end{subfigure}
\caption{Test phantom reconstruction from spindle interior data, noise level $\epsilon=0.05$, 2000 iterations.}
\label{F8}
\end{figure}

Figures \ref{F3}--\ref{F8} show test phantom reconstructions from spindle and spindle interior transform data with varying levels of added noise. In the absence of noise (figures \ref{F3} and \ref{F6}) the reconstructions are ideal. However in the presence of noise we notice a harsher degradation in image quality when reconstructing with spindle interior data. This is because the spindle interior problem is inherently more ill posed due to the extra differentiation step required in the inversion process. At a low noise level (figure \ref{F7}), while the shape of the objects can be deciphered and the ball appears to be hollow, the reconstructions are not clear, in particular the lower density block and the metal sheet start to become lost in the reconstruction. At a higher noise level (figure \ref{F8}), the artefacts in the reconstruction are severe, the ball no longer appears to be hollow and the metal sheet fails to reconstruct. So, although the practical advantages of using a polychromatic source such as an X-ray tube or linear accelerator are clear (i.e. reduced data acquisition time, easy and safe to use etc.), a low noise level would need to be maintained to achieve a satisfactory image quality. The spindle transform inversion performs better when there is added noise. At a low noise level (figure \ref{F4}), the size and shape of the objects is clear and the image contrast is good. At a higher noise level, the background noise in the image is amplified and the ends of metal sheet are hard to identify. We notice, in the presence of noise, that the thinner object (namely the metal sheet) is hardest to reconstruct. This is as we'd expect since smaller densities are determined by the higher frequency harmonic components which degrade faster with noise.

\begin{figure}[h]
\hspace*{-2cm}
\begin{subfigure}{0.5\textwidth}
\includegraphics[width=1.2\linewidth, height=10.5cm]{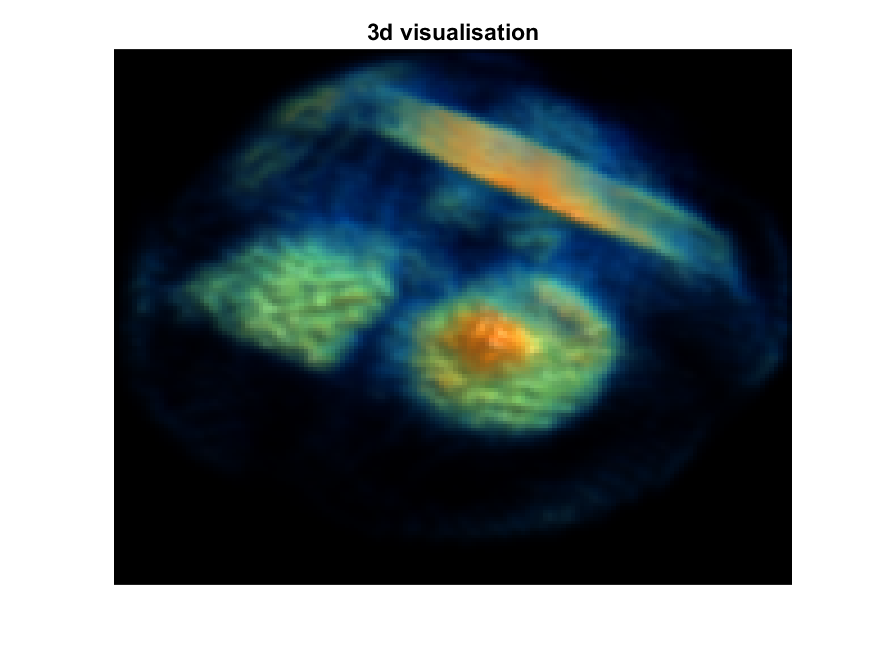} 
\end{subfigure}\hspace{15mm}
\begin{subfigure}{0.5\textwidth}
\includegraphics[width=1.2\linewidth, height=10.5cm]{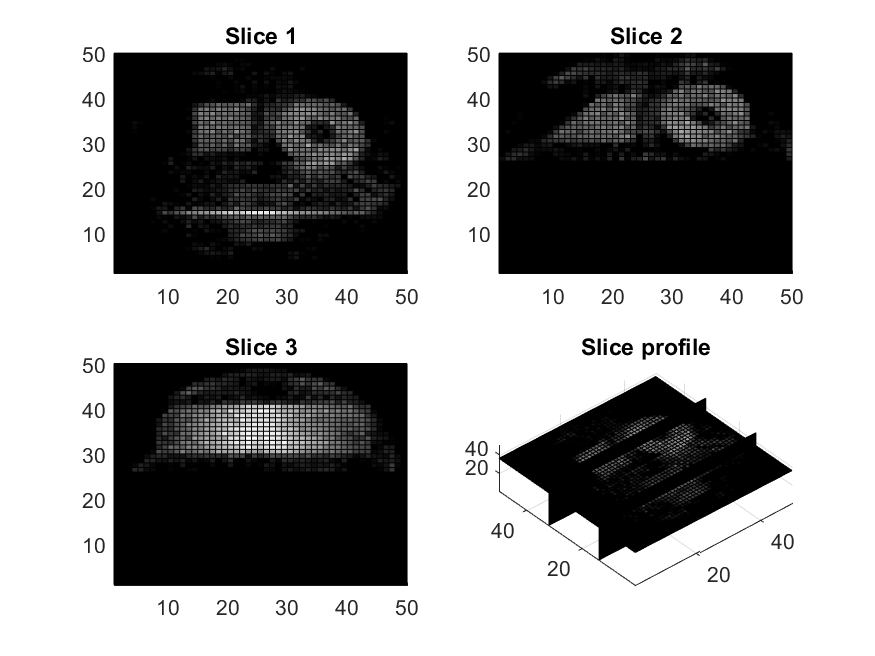}
\end{subfigure}
\caption{Small, low effective $Z$ test phantom reconstruction from spindle transform data with attenuation effects added and noise level $\epsilon=0.01$, 2000  iterations.}
\label{F9}
\hspace*{-2cm}
\begin{subfigure}{0.5\textwidth}
\includegraphics[width=1.2\linewidth, height=10.5cm]{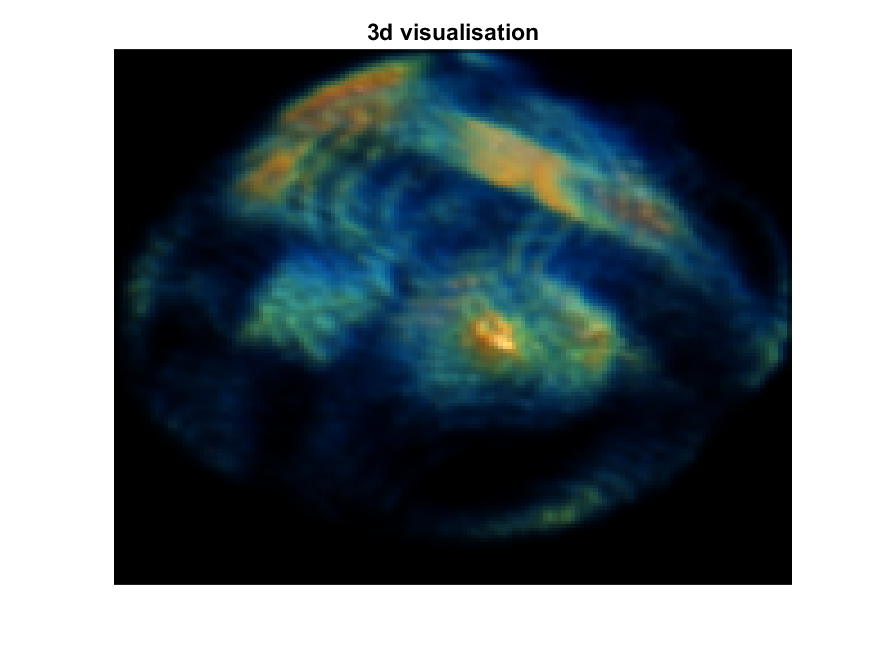} 
\end{subfigure}\hspace{15mm}
\begin{subfigure}{0.5\textwidth}
\includegraphics[width=1.2\linewidth, height=10.5cm]{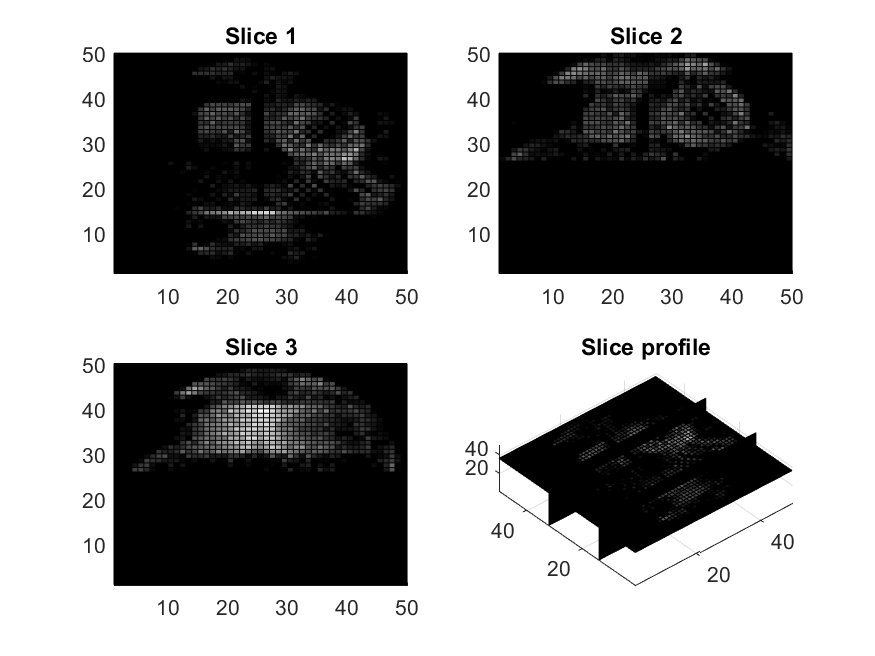}
\end{subfigure}
\caption{Large, high effective $Z$ test phantom reconstruction from spindle transform data with attenuation effects added and noise level $\epsilon=0.01$, 2000 iterations.}
\label{F10}
\end{figure}

To simulate data for the reconstructions presented in figures \ref{F1}--\ref{F8} we applied the discrete operator $A$ to the vector of test phantom pixels $x$ and added Gaussian noise. We now include the effects due the attenuation of the incoming and scattered rays in our simulated data and see how this effects the quality of our reconstruction. For our first example, we set the size of the scanning cube to be $20\times 20\times 20\text{cm}^3$ (each pixel is $4\times 4\times 4$mm$^3$) and the phantom materials are polyethylene (low density block), water (the stairs), rubber (the ball) and Aluminium (the metal sheet). We model the source as Co60, a monochromatic gamma ray source widely used in security screening applications (e.g. screening of freight shipping containers), which has an emission energy of 2824keV. We also add 1$\%$ Gaussian noise as in equation (\ref{error}) to simulate random error as well as the systematic error due to attenuation effects. Our results are presented in figure \ref{F9}. Here the phantom reconstruction is satisfactory with only minor artefacts appearing in the image. So if we scan a low effective $Z$ target of a small enough size with a high energy source, the effects due to attenuation can be neglected while maintaining a satisfactory image quality. For our next example, we set the scanning cube size to be $50\times 50\times 50\text{cm}^3$ (each pixel is $1\times 1\times 1$cm$^3$) and the phantom materials are Teflon (low density block), PVC (stairs), polyoxymethylene (ball) and steel (the metal sheet). The source is Co60 as in our last example and again we add a further 1$\%$ Gaussian noise to the simulated data after we have accounted for the attenuative effects. Our results are presented in figure \ref{F10}. Here, with a larger, higher density target, the effects due to attenuation are significant and the artefacts in the reconstruction are more severe.

\clearpage

\clearpage
\section{Conclusions and further work}
We have presented a new acquisition geometry for three dimensional density reconstruction in Compton imaging with a monochromatic source and introduced a new spindle transform and a generalization of the spindle transform for the surfaces of revolution of a class of symmetric $C^1$ curves. The generalized spindle transform was shown to be injective on the domain of smooth functions $f$ supported on a the intersection of a hollow ball with the upper half space $x_3>0$. In section \ref{spindlesec} it was shown that our problem could be decomposed into a set of one dimensional inverse problems to solve for the harmonic coefficients of a given density, which we then went on to solve via the explicit inversion of a class of Volterra integral operators. Later in section \ref{TI} we considered the problem for a polychromatic source and introduced a new spindle interior transform and proved its injectivity on the set of smooth functions compactly supported on the intersection of a hollow ball and $x_3>0$. We also considered the exterior problem for backscattered photons with a monochromatic and polychromatic photon source, where in section \ref{applesection} we introduced a new apple and apple interior transform. Their injectivity was proven on the set of smooth functions compactly supported on the intersection of the exterior of the unit ball and $x_3>0$. We note that in Palamodov's paper on generalized Funk transforms \cite{pal2}, although he provides an explicit inverse for a fairly general family of integral transforms over surfaces in three dimensional space, the spindle transform is excluded from this family of transforms.

In section \ref{physics} we discussed a possible approach to the physical modelling of our forward problem for both a monochromatic and polychromatic source. In the monochromatic source case, we found that the Compton scattered intensity resembled a weighted spindle transform (as in Theorem \ref{maintheorem}) which could be solved explicitly via repeated approximations. In the polychromatic source case, with a more accurate forward model, we found that an analytic reconstruction was not possible. So we suggested a simplified model in order to obtain an analytic expression for the density, and gave some simple error estimates for our approximation.

Test phantom reconstructions were presented in section \ref{results} using simulated datasets from spindle and spindle interior transform data with varying levels of added pseudo random Gaussian noise. When reconstructing with spindle transform data the reconstructions were satifactory for noise levels up to $5\%$. We saw a harsher reduction in image quality in the presence of added noise when reconstructing from spindle interior data. This was as expected given the increased instability of the spindle interior problem.

In the latter part of section \ref{results}, we provided further reconstructions of our test phantom image in the presence of a systematic error due to the attenuative effects of the incoming and scattered rays. We found, when scanning a small low effective $Z$ target with a high energy source, that the attenuative effects could be neglected while maintaining a good image quality. We also gave an example where the target materials were larger and higher effective $Z$. Here the effects due to attenuation were significant and we saw a more severe reduction in the image quality.

As the stability of the spindle transform has not yet been addressed, in future work we aim to analyse the spindle transform from a microlocal perspective to investigate whether this can shed some light on its stability. Although the injectivity of the exterior problem is covered here, simulations of an exterior density reconstruction are left for future work.

\section*{Acknowledgements}
The authors would like to thank Rapiscan and the EPSRC for the CASE award funding this project, and WL would also like to thank the Royal Society for the Wolfson Research Merit Award that also contributed to this project.

\end{document}